\documentclass[12pt, dvips, twoside]{amsart}
\usepackage{amssymb,latexsym,bbm,graphicx,epsfig,epic,eepic,oldgerm}
\usepackage{psfrag}

\theoremstyle{plain}
\newtheorem{theorem}{Theorem}[section]
\newtheorem{lemma}[theorem]{Lemma}
\newtheorem{proposition}[theorem]{Proposition}
\newtheorem{corollary}[theorem]{Corollary}
\newtheorem{remark}[theorem]{Remark}
\newtheorem{remarks}[theorem]{Remarks}

\newtheorem{example}[theorem]{Example}

\newtheorem{definition}[theorem]{Definition}

\newcommand{\proofend}{\hspace*{\fill} $\Box$\\}
\newcommand{\diam}{\hspace*{\fill} $\Diamond$\\}

\def\s{\smallskip}
\def\m{\medskip}
\def\eps{\varepsilon}
\def\Ker{\operatorname{ker}}
\def\Im{\operatorname {im}}

\def\Symp{\operatorname{Symp}}
\def\Sympcc{\operatorname{Symp_c}}

\def\Int{\operatorname{Int}}
\def\Ham{\operatorname{Ham}}

\def\length{\operatorname{length}}
\def\supp{\operatorname{supp}}
\def\idd{\operatorname{id}}
\def\ind{\operatorname{ind}}
\def\index{\operatorname{index}}
\def\Crit{\operatorname{Crit}}

\def\End{\operatorname{End}}

\def\Fix{\operatorname{Fix^{\circ}}}

\def\Diam{\operatorname{diam}}
\def\HZ{\operatorname{HZ}}
\def\Hess{\operatorname{Hess}}
\def\genus{\operatorname{genus}}
\def\reg{\operatorname{reg}}
\def\tot{\operatorname{tot}}

\def\ga{\alpha}
\def\gb{\beta}
\def\gg{\gamma}
\def\gd{\delta}
\def\eps{\epsilon}

\def\gf{\varphi}

\def\gl{\lambda}
\def\go{\omega}
\def\gs{\sigma}
\def\gt{\vartheta}

\def\ca{{\mathcal A}}
\def\cb{{\mathcal B}}

\def\ce{{\mathcal E}}
\def\cf{{\mathcal F}}
\def\ch{{\mathcal H}}
\def\cj{{\mathcal J}}
\def\cl{{\mathcal L}}
\def\cm{{\mathcal M}}
\def\cn{{\mathcal N}}

\def\cp{{\mathcal P}}

\def\CC{\mathbbm{C}}

\def\NN{\mathbbm{N}}

\def\RR{\mathbbm{R}}

\def\ZZ{\mathbbm{Z}}

\def\pp{\partial}

\def\ra{\rightarrow}

\def\ni{\noindent}
\def\b{\bigskip}
\def\m{\medskip}

\def\id{\mbox{id}}

\def\proof{\noindent {\it Proof. \;}}

\begin{document}

\begin{titlepage}
\title{Hamiltonian dynamics on convex symplectic manifolds}
  
$   $ \\
$   $ \\

\author{Urs Frauenfelder$^1$}
\address{(U.\ Frauenfelder) ETH Z\"urich, CH-8092 Z\"urich, Switzerland}
\email{ufrauenf@math.ethz.ch}
\author{Felix Schlenk$^2$}
\address{(F.\ Schlenk) ETH Z\"urich, CH-8092 Z\"urich, Switzerland}
\email{schlenk@math.ethz.ch}

\date{\today}

\thanks{$^1$ Partially supported by the Swiss National Foundation}
\thanks{$^2$ Supported by the Swiss National Foundation and the von
Roll Research Foundation}

\end{titlepage}

\begin{abstract}
We study the dynamics of Hamiltonian
diffeomorphisms on convex symplectic manifolds.
To this end we first establish the Piunikhin--Salamon--Schwarz isomorphism
between the Floer homology and the Morse homology of such a manifold, 
and then use this isomorphism to construct a biinvariant metric on the
group of compactly supported Hamiltonian diffeomorphisms analogous to
the metrics constructed by Viterbo and Schwarz.
These tools are then applied to prove and reprove results in
Hamiltonian dynamics.
Our applications comprise a uniform lower estimate for the
slow entropy of a compactly supported Hamiltonian diffeomorphism,
the existence of infinitely many nontrivial periodic points of  
a compactly supported Hamiltonian diffeomorphism of a subcritical Stein
manifold, 
old and new cases of the Weinstein conjecture,
and, most noteworthy,
new existence results for closed orbits of a charge in a magnetic field 
on almost all small energy levels.
We shall also obtain some old and new Lagrangian intersection results.
Applications to Hofer's geometry on the group of compactly supported
Hamiltonian diffeomorphisms will be given in \cite{FS2}.
\end{abstract}

\maketitle

 \markboth{{\rm }}{{}} 

\setcounter{tocdepth}{1}
\tableofcontents

\section{Introduction and main results}

\ni
Consider a $2n$-dimensional compact symplectic manifold 
$(M,\go)$ with non-empty boundary $\pp M$.
The boundary $\pp M$ is said to be {\it convex}\, if 
there exists a Liouville vector field $X$ 
(i.e., $\cl_X \go = d \iota_X \go = \go$) 
which is defined near $\pp M$ and is everywhere transverse to $\pp M$, 
pointing outwards;
equivalently, there exists a $1$-form $\ga$ on $\pp M$ such that $d \ga
= \go |_{\pp M}$ and such that $\ga \wedge (d \ga)^{n-1}$ is a volume
form inducing the boundary orientation of $\pp M \subset M$.

\b
\ni
{\bf Definition {\rm (cf.\ \cite{EG})}.}
(i)
A compact symplectic manifold $(M, \go)$ is {\it convex}\, if it has
non-empty convex boundary.

\s
(ii)
A non-compact symplectic manifold $(M, \go)$ is {\it convex}\, if 
there exists an increasing sequence of compact convex submanifolds 
$M_i \subset M$ exhausting $M$, that is,
\[
M_1 \subset M_2 \subset \dots \subset M_i \subset \dots \subset M
\quad \text{ and } \quad
\bigcup_i M_i = M . 
\]  
A symplectic manifold $(M, \go)$ is {\it exact}\, if $\go = d\gl$
and {\it weakly exact}\, if $[\go]$ vanishes on $\pi_2(M)$.

\b
\ni
{\bf Examples.}
%
{\bf 1. Cotangent bundles.}
Recall that every cotangent bundle $T^*N$ over a smooth manifold $N$ carries a
canonical symplectic form $\go_0 = -d \gl$, where $\gl = \sum p_i dq_i$ in
canonical coordinates $(q,p)$. 
The $R$-disc bundles 
\[
T_R^*N \,=\, \left\{ (q,p) \in T^*N \mid \left| p \right| \le R \right\}
\]
over a closed Riemannian manifold $N$ 
and $T^*N = \bigcup_{k \in \NN} T_k^*N$ 
are examples of exact convex symplectic manifolds.
A larger class of examples are

\m
\ni
{\bf 2. Stein manifolds.}
A {\it Stein manifold}\, is a triple $(V,J,f)$ where
$(V,J)$ is an open complex manifold and $f \colon V \ra \RR$ 
is a smooth function which is exhausting and $J$-convex.
``Exhausting'' means that $f$ is bounded from below and proper, and 
``$J$-convex'' means that the $2$-form
\[
\go_f \,=\, -d \left( d f \circ J \right)
\]
is a $J$-positive symplectic form, i.e.,
$\go_f (v,Jv) >0$ for all \text{$v \in TV \setminus \{ 0 \}$}.
We denote by $g_f (\cdot, \cdot) = \go_f (\cdot, J \cdot)$ the induced
K\"ahler metric on $V$, and by $X_f$ the gradient vector filed of $f$ with
respect to $g_f$.
We {\it do not}\, assume that $X_f$ is complete;
in particular, $(V, \go_f)$ can have finite volume.
In any case, 
\begin{equation}  \label{id:Lie}
\cl_{X_f} \go_f \,=\, d \iota_{X_f} \go_f \,=\, 
- d  \left( g_f (X_f, J \cdot) \right) \,=\, 
-d \left( d f \circ J \right) \,=\, \go_f .
\end{equation}
A {\it Stein domain}\, in $(V,J,f)$ is a subset 
$V_R = \left\{ x \in V \mid f (x) \le R \right\}$ 
for a regular value $R \in \RR$.
In view of \eqref{id:Lie}, every Stein domain is an exact compact convex 
symplectic manifold, and so every Stein manifold is an exact convex
symplectic manifold.
We refer the reader to \cite{E1,E2,EG} for foundations of the symplectic
theory of Stein manifolds.

\m
\ni
{\bf 3.}
(i) 
Let $N$ be a closed oriented surface equipped with a Riemannian metric
of constant curvature $-1$, and let $\gs$ be the area form on $N$.
We endow the cotangent bundle $\pi \colon T^*N \ra N$ with the
twisted symplectic form $\go_{\gs} = \go_0 - \pi^* \gs$.
It is shown in \cite{Mac} that $\go_\gs$ is exact on $M = T^*N \setminus
N$ and that $M$ carries a vector field $X$ such that $\cl_X \go_\gs =
\go_\gs$ and such that $X$ is a Liouville vector field on 
\[
M_i \,=\, \left\{ (q,p) \in T^*N \mid \tfrac 1i \le \left| p \right| \le
i \right\}
\]
whenever $i \ge 2$.
Since $H_3 (M_i) = \ZZ$, the manifolds $M_i$, $i \ge 2$, are exact
compact convex symplectic manifolds which are not Stein domains, and $M
= \bigcup_{i \ge 2} M_i$ is an exact convex symplectic manifold which is
not Stein.
Smoothing the boundaries of $k$-fold products $\times_{k} M_i$, $i \ge
2$, we obtain such examples in dimension $4k$ for all $k \ge 1$. 


\s
(ii) Symplectically blowing up a Stein manifold of dimension at least
$4$ at finitely many points we obtain a convex symplectic manifold which
is not weakly exact.

\m
\ni
{\bf 4.} 
A product of convex symplectic manifolds does not need to be convex.
Let $N$ be a closed orientable surface different from the torus,
and let $\gs$ be a $2$-form on $N$.
As we shall see in Lemma~\ref{l:sigma:convex},
the cotangent bundle $T^*N$ endowed with the symplectic form $\go_\gs =
\go_0 - \pi^* \gs$ is convex.
For homological reasons, the product of $\left( T^*N, \go_\gs \right)$
with the convex symplectic manifold $\left( T^* S^1 , \go_0 \right)$ is,
however, convex only if $\gs$ is exact.
We shall be confronted with such non-convex manifolds in our search for
closed trajectories of magnetic flows on surfaces.
We shall therefore develop our tools for symplectic manifolds which away
from a compact subset look like a product of convex symplectic manifolds.
\diam

Throughout we identify $S^1 = \RR / \ZZ$.
Given any symplectic manifold $(M, \go)$, 
we denote by $\ch_c (M)$ the set of $C^2$-smooth functions $S^1 \times M \ra
\RR$ whose support is compact and contained in
$S^1 \times \left( M \setminus \pp M \right)$.
The Hamiltonian vector field of $H \in \ch_c (M)$ defined by
\[
\go \left( X_{H_t}, \cdot \right) \,=\, d H_t \left( \cdot \right)
\]
generates a flow $\gf_H^t$. The set of time-$1$-maps $\gf_H$ form the
group
\[
\Ham_c (M, \go) \,:=\, \left\{ \gf_H \mid H \in \ch_c (M) \right\}
\]
of $C^1$-smooth compactly supported Hamiltonian diffeomorphisms of $(M, \go)$.
Many of our results will apply to those Hamiltonian diffeomorphisms
whose support can be disjoined from itself.
We thus make the

\b
\ni
{\bf Definition.}
A compact subset $A$ of a symplectic manifold $(M, \go)$ is {\it
displaceable}\, 
if there exists $\gf \in \Ham_c(M,\go)$ such that $\gf (A) \cap A =
\emptyset$.  

\b
\ni
{\bf Example.}
Every compact subset of a symplectic manifold of the form $\left( M
\times \RR^2, \go \times \go_0 \right)$ is displaceable.

\b
Our main tools to study Hamiltonian systems on convex symplectic
manifolds will be the Piunikhin--Salamon--Schwarz isomorphism  
and the Schwarz metric.
Before explaining these tools, we describe their applications.
While some applications recover or generalize well-known results, 
many are new; 
all of them, however, are straightforward consequences of the main tools.
In this introduction we give samples of our applications, and we refer
to Sections~\ref{growth} to \ref{lag} and to the appendix for
stronger results.

\b
\ni
{\bf 1. A lower bound for the slow length growth}

\s
\ni
Consider a weakly exact symplectic manifold $(M, \go)$.
For $H \in \ch_c(M)$ the set of contractible $1$-periodic orbits of
$\gf_H^t$ is denoted by $\cp_H$, and the symplectic action $\ca_H(x)$ of
$x \in \cp_H$ is defined as
\begin{equation}  \label{def:af}
\ca_H (x) \,=\, - \int_{D^2} \bar{x}^* \go - \int_0^1 H(t, x(t)) \,dt 
\end{equation}
where $\bar{x} \colon D^2 \ra M$ is a smooth extension of $x$ to the unit disc.
Since $[\go] |_{\pi_2(M)} =0$, the integral $\int_{D^2} \bar{x}^* \go$ does not
depend on the choice of $\bar{x}$.

\b
\ni
{\bf Theorem 1.}
{\it
Assume that $(M, \go)$ is a weakly exact convex symplectic manifold.
Then for every Hamiltonian function $H \in \ch_c(M)$ generating a
non-identical Hamiltonian diffeomorphism $\gf_H \in \Ham_c(M,\go)$ there
exists $x \in \cp_H$ such that $\ca_H(x) \neq 0$.
}
\b

\ni
Theorem~1 is used in \cite{FS} to give a uniform lower bound for the
slow length growth of Hamiltonian diffeomorphisms of exact convex symplectic
manifolds $(M, d\gl)$.
Fix a Riemannian metric $g$ on such a manifold and denote by $\Sigma$
the set of smooth embeddings $\gs \colon [0,1] \ra M$. 
We define the slow length growth $s (\gf) \in [0, \infty]$ 
of a Hamiltonian diffeomorphism $\gf \in \Ham_c (M, \go)$ by
\[
s (\gf) \,=\, \sup_{\gs \in \Sigma} \liminf_{n \ra \infty} 
        \frac{ \log \length_g \left( \gf^n (\gs) \right) }{ \log n} .
\]
Notice that $s (\gf)$ does not depend on the choice of $g$.
We refer to \cite{FS} for motivations to consider this invariant.
Following an idea of Polterovich, \cite{P1}, we use Theorem~1 in
\cite{FS} to show

\b
\ni
{\bf Corollary 1.}
{\it
Assume that $(M, d\gl)$ is an exact convex symplectic manifold.
Then $s (\gf) \ge 1$ for any 
\text{$\gf \in \Ham_c \left( M, d \gl \right) \setminus \left\{ \idd
\right\}$}.
}

\b
\ni
It in particular follows that the group $\Ham_c \left( M, d \gl \right)$
has no torsion.

\b
\ni
{\bf 2. Infinitely many periodic points of Hamiltonian diffeomorphisms}

\s
\ni
We consider again a weakly exact convex symplectic manifold $(M, \go)$.
A {\it periodic point}\, of $\gf_H \in \Ham_c(M, \go)$ is a point $x \in
M$ such that $\gf_H^k(x) =x$ for some $k \in \NN$.
We say that a periodic point $x$ is {\it trivial}\, if $\gf_H^t (x)=x$
and $H_t(x) = 0$ for all $t \in \RR$. 
Since $H \in \ch_c(M)$, $\gf_H$ has many trivial periodic points.
The {\it support}\, $\supp \gf_H$ of a Hamiltonian diffeomorphism 
$\gf_H$ is defined as $\bigcup_{t \in [0,1]} \supp \gf_H^t$. 
It has been proved by Schwarz, \cite{Sch}, in the context of closed weakly
exact symplectic manifolds that if $\supp \gf_H$ is displaceable, 
then $\gf_H$ has infinitely many  geometrically distinct periodic points.
We shall prove an analogous result in our situation.

\b
\ni
{\bf Theorem 2.}
{\it 
Consider a weakly exact convex symplectic manifold $(M, \go)$.
If the support of 
\text{$\gf_H \in \Ham_c (M, \go) \setminus \{ \idd \}$}
is displaceable,
then $\gf_H$ has infinitely many nontrivial
geometrically distinct periodic points corresponding to contractible 
periodic orbits.  
}

\b
\ni
Theorem~2 covers Proposition~4.13 (2) of \cite{V1} 
stating that {\it any}\, non-iden\-ti\-cal compactly supported 
Hamiltonian diffeomorphisms of $\left( \RR^{2n}, \go_0 \right)$ 
has infinitely many nontrivial geometrically distinct periodic points,
see also Theorem~11 in Chapter~5 of \cite{HZ}.
In fact, this is true for all subcritical Stein manifolds.

\b
\ni
{\bf Example (Subcritical Stein manifolds).}
Let $(V,J,f)$ be a Stein manifold.
If $f \colon V \ra \RR$ is a {\it Morse}-function, then
$\index_x(f) \le \frac{1}{2} \dim_{\RR} V$ for all critical points $x$
of $f$. A Stein manifold $(V,J,f)$ is called {\it subcritical}\, if $f$
is Morse and $\index_x(f) < \frac{1}{2} \dim_{\RR} V$ for all 
critical points $x$. 
The simplest example of a subcritical Stein manifold is $\CC^n$ 
endowed with its standard complex structure $J$ and the 
$J$-convex function $f(z_1, \dots, z_n) = |z_1|^2+ \dots +|z_n|^2$.
\diam

\ni
It has been recently shown by Cieliebak, \cite{Ci}, that every
subcritical Stein manifold is symplectomorphic to the product of 
a Stein manifold with $\left( \RR^2, \go_0 \right)$,
and so every compact subset of a subcritical Stein manifold is
displaceable.
We shall not use this difficult result but will 
combine Theorem~2 with a result from \cite{BC} to conclude

\b
\ni
{\bf Corollary 2.}
{\it Any compactly supported non-identical Hamiltonian diffeomorphism of
a subcritical Stein manifold has infinitely many nontrivial
geometrically distinct periodic points corresponding to contrac\-ti\-ble 
periodic orbits.  
}

\b
\ni
{\bf 3. The Weinstein conjecture}

\s
\ni
Another immediate application of our methods is a proof of the Weinstein
conjecture for a large class of hypersurfaces of contact type.
We recall the

\b
\ni
{\bf Definition.}
A $C^2$-smooth compact hypersurface $S$ without boundary of a symplectic
manifold $(M, \go)$ is called {\it of contact type}\, if there exists a
Liouville vector field $X$ which is defined in a neighbourhood of $S$
and is everywhere transverse to $S$. 
A {\it characteristic}\, on $S$ is an embedded circle in $S$ all of whose
tangent lines belong to the distinguished line bundle
\[
\cl_S \,=\, \left\{ (x, \xi) \in TS \mid \go(\xi, \eta) =0 \text{ for
all } \eta \in T_x S \right\} .
\]

\b
\ni
{\bf Theorem~3.}
{\it
Consider a weakly exact convex symplectic manifold $(M, \go)$,
and let $S \subset M \setminus \pp M$ be a displaceable
$C^2$-smooth hypersurface of contact type.
Then $S$ carries a closed characteristic which is contractible in $M$.
}

\b
\ni
Theorem~3 implies a result first proved by Viterbo, \cite{V2}. 

\b
\ni
{\bf Corollary~3.}
{\it
Any $C^2$-smooth hypersurface of contact type in a subcritical Stein
manifold $(V,J,f)$ carries a closed characteristic which is contractible
in $V$.
}
 
\b
\ni
We shall also obtain new existence results for closed characteristics
nearby a given hypersurface.
Roughly speaking, our methods allow to generalize the results which can
be derived from the Hofer--Zehnder capacity for hypersurfaces
in $\RR^{2n}$ to displaceable hypersurfaces in weakly exact convex
symplectic manifolds;
in addition, the closed characteristics found are contractible, and their
reduced actions are bounded by twice the displacement energy of the supporting
hypersurface.
We refer to Section~\ref{weinstein} for the precise results.

\b
\ni
{\bf 4. Closed trajectories of a charge in a magnetic field}

\s
\ni
Consider a Riemannian manifold $(N,g)$ of dimension at least $2$.
The motion of a unit charge on $(N,g)$ subject
to a magnetic field derived from a potential $A \colon N \ra TN$ can be
described as the Hamiltonian flow of the Hamiltonian 
$(p,q) \mapsto \frac12 \left| p-\ga \right|^2$ on $\left( T^*N, \go_0 \right)$ 
where $\ga$ is the $1$-form $g$-dual to $A$ and 
where again $\go_0 = -d \gl$ and $\gl = \sum_i p_i dq_i$.
The fiberwise shift $(q,p) \mapsto \left( q, p-\ga(q) \right)$
conjugates this Hamiltonian system with the Hamiltonian system
\begin{equation}  \label{e:Hmag}
H \colon \left( T^*N, \go_{\gs} \right) \,\ra\, \RR, \quad\, H(q,p) =
\frac 12 \left| p \right|^2 ,
\end{equation}
where $\gs = d \ga$ and
where the twisted symplectic form $\go_\gs$ is given by 
$\go_\gs = \go_0 - \pi^* \gs = - d \left( \gl + \pi^* \ga \right)$.
The system \eqref{e:Hmag} is a model for various other problems in
classical mechanics and theoretical physics,
see \cite{No, Ko}.

A trajectory of a charge on $(N,g)$ in the magnetic field $\gs$ has
constant speed, and closed trajectories $\gg$ on $N$ of speed $c >0$
correspond to closed orbits of \eqref{e:Hmag}
on the energy level $E_c = \left\{ H = c^2 /2 \right\}$.
An old problem in Hamiltonian mechanics asks for closed orbits 
on a given energy level $E_c$, see \cite{Gi0}.
We denote by $\cp^\circ \left( E_c \right)$ the set of closed
trajectories on $E_c$ which are contractible in $T^*N$;
notice that $\cp^\circ \left( E_c \right)$ is the set of closed orbits
on $E_c$ which project to contractible closed trajectories on $N$, and 
that if $\dim N \ge 3$, the orbits in $\cp^\circ \left( E_c \right)$
are contractible in $E_c$ itself.

\b
\ni
{\bf Theorem~4.A.}
{\it
Consider a closed manifold $N$ endowed with a $C^2$-smooth Riemannian
metric $g$ and an exact $2$-form $\gs$ which does not vanish
identically. 
There exists $d>0$ such that $\cp^\circ \left( E_c \right) \neq
\emptyset$ for almost all $c \in \;]0, d]$.
}

\b
``Almost all'' refers to the Lebesgue measure on $\RR$.
The number $d>0$ has a geometric meaning: If the Euler characteristic
$\chi(N)$ vanishes, $d$ is the supremum of the real numbers $c$ for
which the sublevel set
\[
H^c \,=\, \left\{ (q,p) \in T^*N \mid H(q,p) = \tfrac 12 |p|^2 \le c
\right\} 
\]
is displaceable in $\left( T^*N, \go_\gs \right)$, and if $\chi(N)$ does
not vanish, $d$ is defined via stabilizing \eqref{e:Hmag} by 
$\left( T^* S^1, dx \wedge dy \right) \ra \RR$, $(x,y) \mapsto \frac 12 |y|^2$.
Theorem~4.A 
generalizes a result of Polterovich \cite{P0} and Macarini
\cite{Mac} who proved $\cp^\circ (E_c) \neq \emptyset$ for a sequence
$c \ra 0$.

\m
If the magnetic field on $(N,g)$ cannot be derived from a potential, the
motion of a unit charge in this field is still described by
\eqref{e:Hmag},
where now $\gs$ is a closed but not exact $2$-form on $N$, see \cite{Gi0}
and again \cite{No, Ko} for further significance of such Hamiltonian systems.
In this introduction we only consider the case that $N$ is
$2$-dimensional. Since $H^2(N;\RR) =0$ if $N$ is not orientable, 
we can assume that $N$ is orientable.

\b
\ni
{\bf Theorem~4.B.}
{\it
Assume that $N$ is a closed orientable surface endowed with a $C^2$-smooth
Riemannian metric $g$ and a closed $2$-form $\gs \neq 0$.
\begin{itemize}
\item[(i)]
If $N$ is a $2$-sphere, there exists $d>0$ such that $\cp^\circ
\left(E_c\right) \neq \emptyset$ for a dense set of values $c \in \;]0,d]$.
\item[(ii)]
If $\genus (N) \ge 2$, there exists $d>0$ such that 
$\cp^\circ \left(E_c\right) \neq \emptyset$ for almost all $c \in \;]0,d]$.
\end{itemize}
}

\m
\ni
Theorem~4.B is new in case that $\gs$ is not symplectic.
We refer to Section~\ref{ss:magnetic:4B}
for a result containing Theorems~4.A and 4.B as special cases and   
to Section~\ref{ss:state} for a comparison of
ours with previous existence results for closed trajectories of a charge in a
magnetic field.

\b
\ni
{\bf 5. Lagrangian intersections}

\s
\ni
Our methods will provide a concise proof of a Lagrangian
intersection result covering some well known as well as some new cases.

\b
\ni
{\bf Theorem~5.}
{\it
Consider a weakly exact convex symplectic manifold $(M, \go)$, and let
$L \subset M \setminus \pp M$ be a closed Lagrangian submanifold such that
\begin{itemize}
\item[(i)]
the injection $L \subset M$ induces an injection $\pi_1(L) \subset
\pi_1(M)$;  
\s
\item[(ii)]
$L$ admits a Riemannian metric none of whose closed geodesics is contractible.
\end{itemize}
Then $L$ is not displaceable.
} 

\newpage
\ni
{\bf The Schwarz metric}

\s
\ni
We shall derive the above results from a biinvariant spectral metric on
the group $\Ham_c(M, \go)$ of compactly supported Hamiltonian
diffeomorphisms of a weakly exact compact convex symplectic manifold 
$(M,\go)$.
We recall that a symplectomorphism $\gt$ of $(M, \go)$ is a
diffeomorphism of $M$ such that $\gt^* \go = \go$.
We denote by $\Symp_c (M,\go)$ the group of symplectomorphisms of
$(M,\go)$ whose support lies in $M \setminus \pp M$.
We also recall that for any symplectic manifold $(M, \go)$,
Hofer's biinvariant metric $d_H$ on $\Ham_c (M, \go)$ is defined by
\[
d_H \left( \gf, \psi \right) = d_H \left( \gf \psi^{-1}, \idd \right), 
\quad\,
d_H \left( \gf, \idd \right) = \inf \left\{ \| H \| \mid \gf = \gf_H \right\},
\]
where
\[
\| H \| \,=\, \int_0^1 
\left( \sup_{x \in M} H(x,t) - \inf_{x \in M} H(x,t) \right) dt .
\]
It is shown in \cite{LM} that $d_H$ is indeed a metric.

\b
\ni
{\bf Theorem~7.}
{\it
Assume that $(M, \go)$ is a weakly exact compact convex symplectic
manifold.
There exists a function $\gg \colon \Ham_c(M, \go) \ra [0, \infty[$ such
that
\begin{itemize}
\item[(i)]
$\gg (\gf) =0$ if and only if $\gf = \idd$;
\item[(ii)]
$\gg (\gf \psi) \le \gg (\gf) + \gg (\psi)$;
\item[(iii)]
$\gg ( \gt \gf \gt^{-1}) = \gg (\gf)$ for all $\gt \in \Symp_c (M, \go)$;
\item[(iv)]
$\gg (\gf) = \gg \left( \gf^{-1} \right)$;
\item[(v)]
$\gg (\gf) \le d_H \left( \gf, \idd \right)$.
\end{itemize}
}

\b
\ni
In other words, $\gg$ is a symmetric invariant norm on $\Ham_c(M, \go)$.
The {\it Schwarz metric}\, $d_S$ defined by
\[
d_S (\gf, \psi) \,=\, \gg \left( \gf \psi^{-1} \right) 
\]
is thus a biinvariant metric on $\Ham_c(M, \go)$ such that $d_S \le d_H$.
While the Hofer metric is a Finsler metric, the Schwarz metric is a
spectral metric in the sense that $\gg \left( \gf \right)$ is the
difference of two action values of $\gf$.
This property and the property that $\gg \left( \gf_H \right) \le 2\,
\gg \left( \psi \right)$ if $\psi$ displaces the support of $\gf_H$ are
crucial for our applications.
Biinvariant metrics on $\Ham_c(M,\go)$ with these properties 
have been constructed for $\left( \RR^{2n}, \go_0 \right)$
and for cotangent bundles over closed bases by Viterbo \cite{V1} and for 
closed symplectic manifolds by Schwarz \cite{Sch} and Oh \cite{Oh2}.
We shall compare $d_S$ with Viterbo's and Hofer's metric in \cite{FS2}.
There, we shall also use the tools of this paper to study Hofer's
geometry on $\Ham_c(M, \go)$.

\m
The main ingredient in the construction of the Schwarz metric is the 
Piunikhin--Salamon--Schwarz isomorphism (PSS isomorphism, for short)
between the Floer homology and the Morse homology of a weakly exact
compact convex symplectic manifold. 
Floer homology for weakly exact {\it closed}\, symplectic manifolds $(M,\go)$
has been defined in Floer's seminal work \cite{F2,F3,F1,F4}.
It is already shown there that the Floer homology of $(M, \go)$ is
isomorphic to the Morse homology of $M$ and thus to the ordinary homology
of $M$ by considering time independent Hamiltonian functions.
An alternative construction of this isomorphisms was described in
\cite{PSS}; it goes under the name PSS isomorphism.
In the following three sections we establish the PSS isomorphism for
weakly exact compact convex symplectic manifolds $(M, \go)$.
In Sections~\ref{selector} to \ref{metric} we follow \cite{Sch} 
and use our PSS isomorphism to construct the
Schwarz metric $d_S$ on the group $\Ham_c(M, \go)$.
In Section~\ref{ece} we show that the $\pi_1$-sensitive Hofer--Zehnder
capacity is bounded from above by twice the displacement energy.
The last five sections contain our applications.
In the appendix our tools and their applications are extended to 
all convex symplectic manifolds $(M, \go)$ for which the first Chern
class $c_1(\go)$ vanishes on $\pi_2(M)$.

\b
\ni
{\bf Acknowledgements.}
We cordially thank Viktor Ginzburg for introducing us to the
symplectic geometry of magnetic flows.
We also thank Yuri Chekanov, Kai Cieliebak, Urs Lang, 
Leonid Polterovich, Dietmar Salamon, Matthias Schwarz,
Kris Wysocky and Edi Zehnder
for valuable discussions.
This work was done during the second authors stay at FIM of ETH Z\"urich in
the winter term 2002/2003. 
He wishes to thank FIM for its kind hospitality. 

\section{Convexity}  \label{convex}

\ni
We consider a weakly exact compact convex $2n$-dimensional 
symplectic manifold $(M, \go)$.
Choose a smooth vector field $X$ on $M$ which points outwards along $\pp M$
and is such that $\cl_X \go = d \iota_X \go = \go$ near $\pp M$.
For the $1$-form $\ga := (\iota_X \go) |_{\pp M}$ we then have $d \ga =
\go |_{\pp M}$ and $\ga \wedge (d \ga)^{n-1}$ is a volume form inducing
the boundary orientation of $\pp M$. 
Using $X$ we can symplectically identify a neighbourhood of $\pp M$ with
\[
\left( \pp M \times (-2 \eps, 0], d \left( e^r \ga \right) \right)
\]
for some $\eps >0$.
Here, we used coordinates $(x,r)$ on $\pp M \times (-2\eps, 0]$, and in
these coordinates,
$X(x,r) = \frac{\pp}{\pp r}$ on $\pp M \times (-2\eps, 0]$.
We can thus view $M$ as a compact subset of the non-compact symplectic
manifold 
$(\widehat{M},\widehat{\omega})$ defined as
\begin{eqnarray*}
\widehat{M}      &=& M \cup_{\pp M \times \{0\}} \pp M \times [0,\infty), \\
\widehat{\omega} &=& 
   \left\{ \begin{array}{lll} 
          \omega & \text{on} & M, \\
          d \left( e^r \ga \right) & \text{on} & \pp M \times (-2\eps,\infty),
           \end{array}
   \right.
\end{eqnarray*}
and $X$ smoothly extends to $\widehat{M}$ by
\[
\widehat{X} (x,r) := \frac{\partial}{\partial r}, \quad\, 
(x,r) \in \partial M \times (-2\eps, \infty) .
\]
We denote the open ``proboscis'' $\pp M \times (-\eps, \infty)$ by $P_\eps$.
Let $\gf_t$ be the flow of $\widehat{X}$. 
Then $\gf_r(x,0) = (x,r)$ for $(x,r) \in P_\eps$.
We recall that an almost complex structure $\widehat{J}$ on $\widehat{M}$
is called {\it $\widehat{\go}$-compatible}\, if 
\[
\langle \cdot, \cdot \rangle \,\equiv\, 
g_{\widehat{J}} (\cdot, \cdot) \,:=\, 
\widehat{\go} \big( \cdot, \widehat{J} \cdot \big) 
\]
defines a Riemannian metric on $\widehat{M}$.
Following \cite{BPS} we choose an $\widehat{\go}$-compatible 
almost complex structure $\widehat{J}$ on $\widehat{M}$ such that
\begin{eqnarray}  
\widehat{\go} \left( \widehat{X}(x),\widehat{J} (x) v \right) = 0, \;\;
                                   && x \in \pp M, \,\, v \in T_x \pp M,
                                        \label{J3} \\ 
\widehat{\go} \left( \widehat{X}(x), \widehat{J}(x) \widehat{X}(x)
                                   \right) = 1,  \;\; && x \in \pp M, 
                                        \label{J2} \\ 
d \gf_r (x) \widehat{J}(x) = \widehat{J} (x,r) d \gf_r(x), \;\;
             && (x,r) \in P_\eps, \label{J1} 
\end{eqnarray}
We define $f \in C^\infty \left( P_\eps \right)$ by
\begin{equation}\label{f}
f(x,r) := e^r, \quad\, (x,r) \in P_\eps .
\end{equation}
Since $\cl_{\widehat{X}} \widehat{\go} = \widehat{\go}$ on $\pp M \times
(-2\eps, \infty)$, we have $\gf_r^* \widehat{\go} = e^r \widehat{\go}$
on $P_\eps$ for all $r > -\eps$.
This, \eqref{J2} and \eqref{J1} imply that 
\begin{equation}  \label{e:Xf}
\left\langle \widehat{X}(p), \widehat{X}(p) \right\rangle = f(p),   
                                            \quad \, p \in P_\eps .
\end{equation}
Together with \eqref{J3} this implies that
\begin{equation}  \label{e:fX}
\nabla f(p) = \widehat{X} (p), \quad \, p \in P_\eps ,
\end{equation}
where $\nabla$ is the gradient with respect to the metric
$\langle \cdot,\cdot \rangle$. We shall need the following theorem of Viterbo,
\cite{V1}.
\begin{theorem}  \label{t:convex}
For $h \in C^\infty (\RR)$ define
$H \in C^\infty \left( P_\eps\right)$ by
\[
H(p) = h (f(p)), \quad\, p \in P_\eps.
\]
Let $\Omega$ be a domain in $\CC$ and let 
$\widehat{J} \in \Gamma \big( \widehat{M} \times \Omega, \mathrm{End} \big(
T \widehat{M} \big) \big)$ be a smooth
section such that $\widehat{J}_z := \widehat{J} (\cdot,z)$ is an 
$\widehat{\go}$-compatible almost complex structure satisfying 
\eqref{J3}, \eqref{J2} and \eqref{J1}.
If $u \in C^\infty \left( \Omega,P_\eps \right)$ is a solution of 
Floer's equation 
\begin{equation}  \label{e:floer}
\partial_s u(z) + \widehat{J} (u(z),z) \pp_t u(z) = \nabla H(u(z)),\quad \,
z=s+it \in \Omega,
\end{equation}
then
\begin{equation}  \label{viterbo}
\Delta ( f(u) ) = \langle \pp_s u, \pp_s u \rangle
          +h''(f(u))\cdot \pp_s ( f(u) ) \cdot f(u) .
\end{equation}
\end{theorem}

\proof
We abbreviate
$
d^c ( f(u) ) := d ( f(u) ) \circ i = 
\pp_t (f(u)) ds - \pp_s (f(u)) dt$. Then
\begin{equation}  \label{e:lap}
- d d^c ( f(u) ) \,=\, \Delta ( f(u) ) \,ds \wedge dt .
\end{equation}
In view of the identities \eqref{e:Xf}, \eqref{e:fX} and \eqref{e:floer}
we can compute
\begin{eqnarray}  
-d^c (f(u))  
&=& - \left( df(u) \pp_t u \right) ds + \left( df(u) \pp_s u \right) dt 
                                                       \label{e:dc} \\ 
&=& - \left( df(u) \big( \widehat{J}(u,z) \pp_t u \big) \right)dt
       - \left( df(u) \big( \widehat{J}(u,z) \pp_s u \big) \right)ds \notag\\ 
& & + \left( df(u) \big( \pp_s u + \widehat{J}(u,z) \pp_t u \big) \right) dt
    + \left( df(u) \big( \widehat{J}(u,z) \pp_su-\pp_t u \big) \right)
                                                              ds \notag\\ 
&=&   \widehat{\go} \,\big( \widehat{X}(u), \pp_t u \big) \,dt 
    + \widehat{\go} \,\big( \widehat{X}(u), \pp_s u \big) \,ds \notag\\ 
& & + \left\langle \nabla f(u), \nabla H(u) \right\rangle dt
    + \big\langle \nabla f(u), \widehat{J}(u,z) \nabla H(u) \big\rangle
                                                           \,ds \notag\\  
&=&u^* \iota_{\widehat{X}} \widehat{\go} + \big\langle \widehat{X}(u), h'(f(u))
                            \widehat{X} (u) \big\rangle \,dt +0 \notag\\
&=&u^* \iota_{\widehat{X}} \widehat{\go} +h'(f(u)) f(u) dt.  \notag
\end{eqnarray}
Using $d \iota_{\widehat{X}} \widehat{\go} = \cl_{\widehat{X}} \widehat{\go} =
\widehat{\go}$ and again \eqref{e:floer}, we find
\begin{eqnarray*}
d u^* \iota_{\widehat{X}} \widehat{\go} \,=\, u^* \widehat{\go} &=& 
 \widehat{\go} \left( \pp_s u , \widehat{J}(u,z) \pp_s u - \widehat{J}(u,z) \nabla H(u) \right)
                                               ds \wedge dt  \\
     &=& \big( \langle \pp_s u, \pp_s u \rangle  -
            dH(u) \pp_s u \big)\, ds \wedge dt \\
     &=& \big( \langle \pp_s u, \pp_s u \rangle - \pp_s (h(f(u))) \big)\, ds
     \wedge dt. 
\end{eqnarray*}
Together with \eqref{e:dc} it follows that
\begin{eqnarray*}
-d d^c ( f(u) ) &=&
\big( \langle \pp_s u, \pp_s u \rangle - \pp_s (h(f(u))) 
   +\pp_s(h'(f(u))f(u)) \big)\, ds \wedge dt\\
&=& \big( \langle \pp_s u, \pp_s u \rangle 
      + h''(f(u))\cdot \partial_s f(u) \cdot f(u) \big)\, ds \wedge dt ,
\end{eqnarray*}
and so Theorem~\ref{t:convex} follows in view of \eqref{e:lap}.
\proofend

\begin{remark}[Time-dependent Hamiltonian]  \label{sconvex}
{\rm
Repeating the calculations in the proof of
Theorem~\ref{t:convex}, one shows the following more general result.} 
{\it Let 
$h \in C^\infty(\mathbb{R}^2,\mathbb{R})$ and define
$H \in C^\infty \left( P_\eps \times \RR \right)$ by
\[ 
H(p,s) = h(f(p),s), \quad\, p \in P_\eps, \,\,  s \in \RR . 
\] 
If $\Omega$ is a domain in $\CC$ and if  
$u \in C^\infty \left( \Omega, P_\eps \right)$ is a solution of the
time-dependent Floer equation 
\begin{equation}\label{sfloer} 
\pp_s u(z)+ \widehat{J}(u(z),z) \pp_t u(z) = \nabla H(u(z),s), \quad\,
z=s+it \in \Omega , 
\end{equation} 
then 
\[ 
\Delta (f(u)) = \langle \partial_s u,\partial_s u \rangle +
\partial^2_1h(f(u),s)\cdot \partial_s f(u) \cdot f(u) 
+\partial_1\partial_2 h(f(u),s)\cdot f(u) . 
\] 
} 
\end{remark}

In the following corollary we continue the notation of
Theorem~\ref{t:convex}.

\begin{corollary}[Maximum Principle]  \label{maximum}
Assume that $u \in C^\infty \left( \Omega, P_\eps \right)$
and that one of the following conditions holds.
\begin{itemize}
 \item[(i)]$u$ is a solution of Floer's equation \eqref{e:floer};
 \item[(ii)]
  $u$ is a solution of the time-dependent Floer equation~\eqref{sfloer}
  and $\partial_1\partial_2 h \geq 0$. 
\end{itemize}
If $f \circ u$ attains its maximum on $\Omega$, then $f \circ u$ is constant.
\end{corollary}

\proof
Assume that $u$ solves \eqref{e:floer}.
We set \[
b(z) \,=\, -h'' (f(u(z))) \cdot f(u(z)) .
\]
The operator $L$ on $C^\infty (\Omega,\RR)$ defined by
$
L(v) = \Delta v + b(z) \pp_s v
$
is uniformly elliptic on relatively compact domains in $\Omega$, and
according to Theorem~\ref{t:convex}, $L(f \circ u) \ge 0$. 
If $f \circ u$ attains its maximum on $\Omega$, 
the strong Maximum Principle \cite[Theorem 3.5]{GT} thus implies that $f
\circ u$ is constant.
The other claim follows similarly from 
Remark~\ref{sconvex} and the second part of \cite[Theorem 3.5]{GT}. 
\proofend

\section{Floer homology}

\ni
The Floer chain complex of a Hamiltonian function is generated by the
$1$-periodic orbits of its Hamiltonian flow,
and the boundary operator is defined by counting perturbed 
pseudo-holomorphic cylinders which converge at both ends 
to generators of the chain complex. 
In the presence of a contact-type boundary the Hamiltonian 
has to be chosen appropriately near the boundary in order to insure 
that the Floer cylinders stay in the interior of the manifold.  

\s
The Reeb vector field $R$ of $\ga$ on $\pp M$ is defined by
\begin{equation}\label{reeb}
\omega_x ( v,R ) = 0 
\,\text{ and }\, 
\omega_x (X,R) = 1 , \quad \,
\, x \in \pp M, \,\, v \in T_x \pp M.
\end{equation}
By \eqref{J3} and \eqref{J2} we have $R = \widehat{J} \,X |_{\pp M}$. This
and \eqref{e:fX} imply that for 
$h \in C^\infty ( \RR )$  the Hamilton equation 
$\dot{x} = X_H (x)$ of $H = h \circ f \colon P_\eps \ra \RR$ 
defined by $\go \left( X_H (x) , \cdot \right) = d H (x)$
restricts on $\pp M$ to
\begin{equation}  \label{hamilton}
\dot{x}(t) \,=\, h'(1) \,R(x(t)) .
\end{equation}
Define $\kappa \in (0, \infty ]$ by
\[ 
\kappa :=
\inf \left\{ c > 0 \mid \dot{x}(t) = c\, R(x(t))
        \text{ has a $1$-periodic orbit} \right\} .
\]
We denote by $\widehat{\ch}$ the set of smooth functions
$\widehat{H} \in C^\infty( S^1 \times \widehat{M} )$ for which there exists
$h \in C^\infty( \RR )$ such that $0 \le h'(\rho) < \kappa$ for all
$\rho \ge 1$ and
$\widehat{H}|_{ S^1 \times \pp M \times [0,\infty) } = h \circ f$;
with this choice of $h$ the restriction of the flow
$\gf_{\widehat{H}}^t$ of $\widehat{H} \in \widehat{\ch}$ to $\pp M
\times [0,\infty)$ has no $1$-periodic solutions.
We introduce the set 
\[
\ch := \left\{ H \in C^\infty ( S^1 \times M ) \mid  
   H = \widehat{H}|_{S^1 \times M} \text{ for some } 
   \widehat{H} \in \widehat{\ch} \right\} 
\] 
of {\it admissible Hamiltonian functions}\, on $M$.
Moreover,
we denote by $\widehat{\cj}$ the set of smooth sections
$\widehat{J} \in \Gamma \big( S^1 \times \widehat{M} ,\mathrm{End}
\big(T\widehat{M} \big) \big)$  
such that for every $t \in S^1$ the section 
$\widehat{J}_t := \widehat{J} (t, \cdot )$ is an $\widehat{\omega}$-compatible
almost complex structure which on  $\pp M \times [0,\infty)$
is independent of the $t$-variable and
satisfies \eqref{J3}, \eqref{J2} and \eqref{J1};
and we introduce the set 
\[
\cj := \left\{ J \in \Gamma \left( S^1 \times M ,\mathrm{End}(TM) \right) \mid 
 J=\widehat{J}|_{S^1 \times M} \text{ for some } \widehat{J} \in \widehat{\cj} \right\} 
\]
of {\it admissible almost complex structures} on $TM$.
A well-known argument shows that the space $\widehat{\cj}$ is connected,
see \cite[Remark 4.1.2]{BPS}. Since the restriction map $\widehat{\cj} \ra \cj$
is continuous, $\cj$ is also connected.

For $H \in \ch$ let $\cp_H$ be the set of contractible
$1$-periodic orbits of the Hamiltonian flow of $H$.
By ``generic'' we shall mean ``belonging to a countable
intersection of sets which are open and dense in the $C^\infty$-topology''. 
For generic $H \in \ch$ for no $x \in \cp_H$
the value $1$ is a Floquet multiplier of $x$, i.e.,
\begin{equation}  \label{e:floer:det}
\det \left( \id - d \gf^1_{H}(x(0) \right) \neq 0.
\end{equation}
Since $M$ is compact, $\cp_H$ is then a finite set.
An admissible $H$ satisfying \eqref{e:floer:det} for all $x \in \cp_H$ 
is called {\it regular}, 
and the set of regular admissible Hamiltonians is denoted
by $\ch_{\reg} \subset \ch$.
For $H \in \ch_{\reg}$ we define 
$CF (M;H)$ to be the $\ZZ_2$-vector space consisting of formal sums 
\[
\xi \,= \sum_{x \in \cp_H} \xi_x\, x, \quad\, \xi_x \in \ZZ_2.
\]
We assume first that the first Chern class 
$c_1 = c_1(\go) \in H^2(M;\ZZ)$ of the bundle $(TM,J)$ vanishes on $\pi_2(M)$.
In this case,
the Conley--Zehnder index $\mu (x)$ of $x \in \cp_H$ is well-defined,
see \cite{SZ}.
We normalize $\mu$ in such a way that for $C^2$-small 
time-independent Hamiltonians,
\[
\mu(x) = 2n - \ind_H(x)
\]
for each critical point $x \in \Crit(H)$; 
here, $\ind_H (x)$ is the Morse index of $H$ at $x$.
The Conley--Zehnder index turns $CF (M;H)$ into the graded
$\ZZ_2$-vector space $CF_*(M;H)$.
For $x,y \in \cp_H$ let $\cm (x,y)$ be the moduli space of Floer 
connecting orbits from $x$ to $y$, i.e., 
$\cm (x,y)$ is the set of solutions $u \in C^\infty( \RR \times S^1, M)$
of the problem
\begin{eqnarray}  \label{prob:floer} 
 \left\{ 
     \begin{array}{c} 
      \partial_s u+J_t(u)(\partial_t u-X_{H_t}(u)) = 0, \\[0.4em]
      \displaystyle\lim_{s \to -\infty}u(s,t)  = x(t), \quad
      \displaystyle\lim_{s \to \infty}u(s,t) = y(t).
    \end{array}
 \right.
\end{eqnarray}
For later use we notice that by a standard computation,
\begin{equation}  \label{est:action0}
\int_{\RR \times S^1} \left| \pp_s u \right|^2  ds dt \,=\, \ca_H(x) -
\ca_H(y) \,\ge\, 0, \
\quad\, u \in \cm (x,y) .
\end{equation}
For generic $J \in \mathcal{J}$ the moduli spaces 
$\cm (x,y)$ are smooth manifolds of dimension $\mu (x) - \mu (y)$
for all $x,y \in \cp_H$, see \cite{SZ}.
Such a $J$ is called {\it $H$-regular}, and a pair $(H,J)$ is called
regular if $H$ is regular and $J$ is $H$-regular. 
The group $\RR$ acts on $\cm (x,y)$ by translation, 
$u(s,t) \mapsto u(s+\tau,t)$ for $\tau \in \RR$. 
Since $[ \go ]$ vanishes on $\pi_2(M)$, there is no bubbling off
of pseudo holomorphic spheres. 
It thus follows from Corollary~\ref{maximum} that if
$\mu(x)-\mu(y)=1$, then the quotient 
$\cm (x,y) / \RR$ is a compact zero-dimensional manifold
and hence a finite set. Set
\[
n(x,y) \,:=\, \# \left\{ \cm (x,y) / \RR \right\} \mod 2.
\]
For $k \in \NN$ we define the Floer boundary operator 
$\pp_k \colon CF_k(M;H) \to CF_{k-1}(M;H)$ as the linear extension of
\[
\pp_k x \,= \sum_{\substack{y \in \cp_H\\\mu(y)=k-1}} n(x,y) \,y
\]
where $x \in \cp_H$ and $\mu(x)=k$. 
Proceeding as in \cite{F3,Sch1} one shows that $\partial^2 = 0$.
The complex $\left( CF_* (M;H), \pp_* \right)$ is called the Floer chain
complex. Its homology
\[
HF_k(M;H,J) \,:=\,  \frac{\Ker \partial_k}{\Im \partial_{k+1}}
\]
is a graded $\ZZ_2$-vector space which does not depend on 
the choice of a regular pair $(H,J)$, see again \cite{F3,Sch1},
and so we can define the Floer homology $HF_*(M)$ by
\[
HF_*(M) \,:=\, HF_*(M;H,J)
\]
for any regular pair $(H,J)$.

\s
In case that $c_1 (\go)$ does not vanish, the moduli spaces
$\cm(x,y)$ for $x,y \in \cp_H$ are still smooth manifolds for generic $J
\in \cj$, but now may contain 
connected components of different dimensions. We denote
by $\cm^1(x,y)$ the union of the $1$-dimensional connected
components of $\cm(x,y)$. Since $[ \omega ]$ vanishes on $\pi_2(M)$,
the space $\cm^1(x,y)/ \RR$ is still compact, and we can define
\[
n(x,y) \,:=\, \#(\cm^1(x,y)/ \RR) \mod 2 .
\]
Proceeding as above we define an ungraded Floer homology whose
chain complex is generated again by the set $\cp_H$ and
whose boundary operator is the linear extension of
\[
\pp x \,=\, \sum_{y \in \cp_H} n(x,y) \,y
\]
where $x \in \cp_H$. 
We shall explain in the appendix how Novikov rings can be used to define
a graded Floer homology even if $c_1$ does not vanish on $\pi_2(M)$.

\subsection*{Products}

As we pointed out in Example~4 of the introduction,
the product of convex manifolds does not need to be convex.
Nevertheless, the Floer homology of a product
of weakly exact compact convex symplectic manifolds 
can still be defined.
In fact, Floer homology can be defined for a yet larger class of compact
symplectic manifolds with corners.
\begin{definition}  \label{def:splitconvex}
{\rm
A compact symplectic manifold with corners $(M,\go)$ is 
{\it split-convex}\, 
if there exist compact convex symplectic manifolds 
$(M_j, \go_j)$, $j = 1, \dots, k$, and a compact subset 
$K \subset M \setminus \pp M$ such that $M = M_1 \times \dots \times
M_k$ and 
\[
\left( M \setminus K, \go \right) \,=\, \left( \left( M_1 \times
\dots \times M_k \right) \setminus K, \go_1 \oplus \dots \oplus \go_k
\right) .
\]
}
\end{definition}

\ni
Consider a weakly exact compact split-convex symplectic manifold
$(M,\go)$, and let $(M_j, \go_j)$, $j = 1, \dots, k$,
be as in Definition~\ref{def:splitconvex}.
For notational convenience, we assume $k=2$.
We specify the set of admissible Hamiltonian functions
$\ch \subset C^\infty \left( S^1 \times M \right)$
and the set of admissible almost complex structures 
$\cj \subset \Gamma \left( S^1 \times M , \End \left( TM \right) \right)$
as follows. 
For $i = 1,2$, let 
$\widehat{M_i} = M_i \cup_{\partial M_i \times \{0\}}\partial M_i
\times [0,\infty)$ be the completion of $M_i$ endowed with the
symplectic form $\widehat{\go}_i$ as in Section~\ref{convex}, and let
$\widehat{\ch_i} \subset C^\infty \big( S^1 \times \widehat{M_i} \big)$ 
and $\widehat{\cj_i} \subset \Gamma \big( S^1 \times \widehat{M_i}, \End
\big( T \widehat{M_i} \big) \big)$ 
be the set of admissible functions and admissible almost complex
structures on $\widehat{M_i}$. 
We define the completion $\big( \widehat{M}, \widehat{\go} \big)$ of
$(M, \go)$ as
\begin{eqnarray*}
\widehat{M}      &=& \widehat{M_1} \times \widehat{M}_2 , \\
\widehat{\omega} &=& 
   \left\{ \begin{array}{lll} 
          \omega & \text{on} & M, \\ [0.2em]
          \widehat{\go}_1 \oplus \widehat{\go}_2 & \text{on} & 
               \big( \widehat{M_1} \times \widehat{M}_2 \big) 
                          \setminus \left( M_1 \times M_2 \right) .
           \end{array}
   \right.
\end{eqnarray*}
We first define the set of admissible functions 
$\widehat{\ch} 
\subset C^\infty \big( S^1 \times \widehat{M} \big)$ as the set of functions 
$\widehat{H} \in C^\infty \big( S^1 \times \widehat{M} \big)$ 
for which there exist 
$\widehat{H_i} \in \widehat{\ch_i}$, $i=1,2$, such that 
\[
\widehat{H}|_{ ( \widehat{M_1} \times \widehat{M_2} ) \setminus
\left( M_1 \times M_2 \right)} \,=\,
 ( \widehat{H_1} + \widehat{H_2} ) \Big|_
 { ( \widehat{M_1} \times \widehat{M_2} ) \setminus
 \left( M_1 \times M_2 \right)} ;
\]
and we then define the set $\ch$ of admissible functions on 
$M$ as the set of functions
$H \in C^\infty \left( S^1 \times M \right)$ for which
there exists $\widehat{H} \in \widehat{\ch}$ such that
\[
H \,=\, \widehat{H}|_M .
\]
Similarly, we first define the set of admissible almost complex structures
$\widehat{\cj}$ as the set of 
$\widehat{J} \in \Gamma \left( S^1 \times
\widehat{M}, \End \big( T \widehat{M} \big) \right)$ for which there exist
admissible almost complex structures $\widehat{J_i} \in
\widehat{\cj_i}$, $i=1,2$, such that
\[
\widehat{J}|_{ ( \widehat{M_1} \times \widehat{M_2} ) \setminus
\left( M_1 \times M_2 \right)} \,=\, 
( \widehat{J_1} \times \widehat{J_2} )
\Big|_{ (\widehat{M_1} \times \widehat{M_2} ) \setminus
\left( M_1 \times M_2 \right)} ;
\]
and we then define the set $\cj$ of admissible almost complex structures on
$M$ as the set of almost complex structures 
$J \in \Gamma \left( S^1 \times M , \End \left( TM \right) \right)$ 
for which there exists $\widehat{J} \in \widehat{\cj}$ such that
\[
J \,=\, \widehat{J}|_M .
\]
Using the maximum principle Corollary~\ref{maximum} factorwise we define
the Floer homology $HF \left( M \right)$
as above. 
If $c_1 \left( \go \right)$ vanishes on $\pi_2 \left( M \right)$, 
then $HF \left( M \right)$ is graded by the Conley--Zehnder
index.

\section{The Piunikhin--Salamon--Schwarz isomorphism}  \label{iso}
 
\ni
We assume again that $(M,\go)$ is a weakly exact compact convex
symplectic manifold.
We first assume that $c_1(\go)$ vanishes on $\pi_2(M)$. 
Let $F \in C^\infty(M)$ be an {\it admissible Morse function}, i.e., 
$F$ is a smooth Morse function for which there exists 
$\widehat{F} \in C^\infty(\widehat{M})$ such that
\begin{equation*}
\widehat{F}|_M = F \,\text{ and }\,
\widehat{F}(x,r)=e^{-r}, \quad \, x \in \partial M,\,\,r \in [0,\infty).
\end{equation*}
The Morse chain complex $CM_*(M;F)$ of $F$ is the $\ZZ_2$-vector space
generated by the critical points of $F$ and graded by the Morse index,
and the boundary operator on $CM_*(M;F)$ is defined by counting 
flow lines of the negative gradient flow of $F$ with respect to a
generic Riemannian metric between critical points of index difference $1$. 
The homology 
\[
HM_*(M) \,=\, HM_*(M;F)
\]
of $CM_*(M;F)$
does not depend on the choice of $F$, cf.\ \cite{Sch1}.
The Piunikhin--Salamon--Schwarz
maps will give us an explicit isomorphism between the Floer homology
$HF_*(M)$ of $(M, \go)$ and the Morse homology $HM_*(M)$ of $M$. 

\s
Choose $H \in \ch_{\reg}$
and an admissible Morse function $F$.
We first construct the Piunikhin--Salamon--Schwarz map 
\[
\phi \colon CM_*(M;F) \,\ra\, CF_*(M;H) .
\]
By definition of $\ch$, there exists $\widehat{H} \in \widehat{\ch}$
such that $H = \widehat{H}|_{S^1 \times M}$ and
 $\widehat{H}|_{S^1 \times \pp M \times [0,\infty)} = h \circ f$
for some $h \in C^\infty(\mathbb{R})$ satisfying $0\leq  h'(1) <\kappa$. 
For $s \in \RR$ choose a smooth family
$h_s \in C^\infty( \RR )$ such that
\begin{itemize}
 \item[(h1)] $h_s=0, \quad s \le 0$,
 \item[(h2)] $\pp_s h'_s \geq 0, \quad s \in \RR$,
 \item[(h3)] $h_s = h, \quad s \ge 1$,
\end{itemize}
and then choose a smooth family 
$\widehat{H}_s \in C^\infty \big( S^1 \times \widehat{M} \big)$ such that
\begin{itemize}
 \item[(H1)] $\widehat{H}_s = 0, \quad s \le 0$,
 \item[(H2)] $\widehat{H}_s|_{S^1 \times \partial M \times [0,\infty)}
  = h_s \circ f, \quad 0\leq s\leq 1$,
 \item[(H3)] $\widehat{H}_s=\widehat{H}, \quad s \ge 1$.
\end{itemize}
We finally define the smooth family $H_s \in C^\infty( S^1 \times M )$ by
\[
H_s \,:=\, \widehat{H}_s|_{S^1 \times M}.
\]
\begin{theorem}  \label{t:4finite}
Let $J^-$ be an $H$-regular admissible almost complex structure
and $J^+$ be an arbitrary admissible almost complex structure.
Consider the space $\cj (J^-,J^+)$ of 
families of admissible almost complex structures
$J_s \in \Gamma( S^1 \times M , \End (TM))$ for which there
exists $s_0=s_0(\widehat{J}_s)>0$ such that $J_s=J^-$ for $s \leq -s_0$ and
$J_s=J^+$ for $s \geq s_0$.
For a generic element $J_s \in \cj (J^-,J^+)$
the moduli space of the problem
\begin{equation}  \label{PSS}
\left\{
\begin{array}{l}
u \in C^\infty(\RR\times S^1,M), \\ [0.4em]
\pp_s u+J_{s,t}(u) \left( \pp_t u-X_{H_{s,t}}(u) \right) =0, \\ [0.4em]
\int_{\RR \times S^1} \left| \pp_s u \right|^2 < \infty ,
\end{array}\right.
\end{equation}
is a smooth finite dimensional manifold. 
Here, $J_{s,t}=J_s(\cdot,t)$ and \text{$H_{s,t}=H_s(\cdot,t)$}.
\end{theorem}

\proof 
We denote the moduli space of solutions of problem~\eqref{PSS}
by $\cm_0$. Choose $\widehat{H}_s \in \widehat{\ch}$
and $\widehat{J}_s \in \widehat{\cj}$ satisfying
$\widehat{H}_s|_{M}=H_s$ and $\widehat{J}_s|_{M}=J_s$.
Instead of $\cm_0$ we first consider the moduli space $\cm$ of solutions 
of the problem
\begin{equation}  \label{PSS2}
\left\{
\begin{array}{l}
u \in C^\infty(\RR \times S^1, \widehat{M}), \\ [0.4em]
\partial_s u+\widehat{J}_{s,t}(u)(\partial_t u-
X_{\widehat{H}_{s,t}}(u))=0, \\ [0.4em]
\int_{\RR \times S^1} \left| \pp_s u \right|^2 < \infty .
\end{array}\right.
\end{equation}
The moduli space $\cm_0$ consists of those
$u \in \cm$ whose image is entirely contained in
$M$. We shall first prove that for generic choice of 
$\widehat{J}_s$, the moduli space $\cm$ is a smooth finite
dimensional manifold. 
We shall then use convexity to prove that
the image of each $u \in \cm$ is entirely contained in $M$ and
hence $\cm_0 = \cm$ is a smooth finite dimensional manifold.

We interpret solutions of \eqref{PSS2} as
the zero set of a smooth section from a Banach manifold $\cb$
to a Banach bundle $\ce$ over $\cb$. To define 
$\cb$ we first introduce certain weighted Sobolev norms. 
Choose a smooth cutoff function $\beta \in C^\infty(\mathbb{R})$
such that $\beta(s)=0$ for $s<0$ and $\beta(s)=1$ for $s>1$. Choose
$\delta>0$ and define $\gamma_\delta \in C^\infty(\mathbb{R})$
by
\[
\gamma_\delta(s) \,:=\, e^{\delta \beta(s) s}.
\]
Let $\Omega$ be a domain in
the cylinder $\RR \times S^1$. 
For $1 \leq p \leq \infty$ and $k \in \NN_0$ 
we define the $\left\| \cdot \right\|_{k,p,\delta}$-norm 
for $v \in W^{k,p}(\Omega)$ by
\[
\left\| v \right\|_{k,p,\delta} \,:= \sum_{i+j \leq k} 
\left\| \gamma_\delta \cdot \pp^i_s\partial^j_t v \right\|_p \,.
\]
We introduce weighted Sobolev spaces
\begin{eqnarray*}
W^{k,p}_\delta(\Omega) &:=& \left\{ v \in W^{k,p}(\Omega) \mid 
\left\| v \right\|_{k,p,\delta} < \infty \right\}  \\
  &=& \left\{ v \in W^{k,p}(\Omega) \mid \gamma_\delta v \in
   W^{k,p}(\Omega) \right\},
\end{eqnarray*}
and we abbreviate
\[
L^p_\delta(\Omega) \,:=\, W^{0,p}_\delta(\Omega).
\]
Let $p>2$ and fix a metric $g$ on $T \widehat{M}$. 
The Banach manifold $\cb=\cb^{1,p}_\delta(\widehat{M})$
consists of $W^{1,p}_\text{loc}$-maps $u$ from the cylinder 
$\RR \times S^1$ to $\widehat{M}$ which satisfy the conditions
\begin{itemize}
 \item[(B1)] 
  There exists a point $m \in \widehat{M}$, 
  a real number $T_1 <0$, and
  $v_1 \in W^{1,p}_\delta \big( (-\infty, T_1) \times S^1, T_p
  \widehat{M} \big)$ such that
  \[
  u(s,t) \,=\, \exp_m (v_1 (s,t)), \quad\, s < T_1 .
  \]
 \item[(B2)]
  There exists $x \in \cp_H \subset C^\infty(S^1,\widehat{M})$,
  a real number $T_2 >0$, and $v_2 \in W^{1,p}_\gd \big( (T_2, \infty)
  \times S^1,x^*TM \big)$
  such that
  \[
  u(s,t) \,=\, \exp_{x(t)} (v_2 (s,t)), \quad\, s > T_2 .
  \]
\end{itemize}
Here, the exponential map is taken with respect to $g_J$. 
Since $\widehat{M}$ has no boundary, $\cb$ is a Banach manifold
without boundary. 
Note that every solution of \eqref{PSS2} lies in $\cb$. 
Indeed, the finite energy assumption in \eqref{PSS2}
guarantees that solutions of \eqref{PSS2} converge exponentially fast at
both ends, see \cite[Section 3.7]{F}.
Let $\ce$ be the Banach bundle over $\cb$ whose fiber over 
$u \in \cb$ is given by
\[
\ce_u \,:=\, L^p_\delta (u^*TM).
\]
We choose $\widehat{J}^-, \widehat{J}^+ \in \widehat{\cj}$ such that
$J^-=\widehat{J}^-|_M$ and $J^+=\widehat{J}^+|_M$. For each smooth family
$\widehat{J}_s$ for which there exists an $s_0>0$ such that
$\widehat{J}_s=\widehat{J}^-$ for $s \leq -s_0$ and
$\widehat{J}_s=\widehat{J}^+$ for $s \geq s_0$
we define the section
$\cf = \cf_{J_s} \colon \cb \to \ce$ by
\[
\cf (u) \,:=\, \pp_s u+\widehat{J}_{s,t}(u)(\pp_t u-X_{\widehat{H}_{s,t}}(u)) .
\]
If $\delta$ is chosen small enough, then the
vertical differential $D \cf$ is a Fredholm operator, see
for example \cite[Section 4.3]{F}. 
One can prove that for generic choice of $J_s$ 
the section $\cf_{J_s}$ intersects the zero section transversally, see
\cite[Section 5]{FHS} and \cite[Section 4.5]{F}.
Hence,
\[
\cm \equiv \cm_{J_s} \,:=\, \cf_{J_s}^{-1}(0)
\]
is a smooth finite dimensional manifold for generic $J_s$. 

It remains to show that $\cm = \cm_0$, 
i.e., the image of every $u \in \cm$ is contained in $M$. 
We first claim that $m: = \lim_{s \to -\infty}u(s,t) \in M$. 
To see this, assume that $m \in \widehat{M} \setminus M$. 
Define $v \colon \CC \to \widehat{M}$ by the conditions
\[
v \left( e^{2 \pi(s+it)} \right) = u(s,t), \quad\,  v(0) = m .
\]
Since every admissible almost complex structure $J$ restricted to
$\widehat{M}\setminus M$ is independent of the $t$-variable, 
$v$ is a pseudo holomorphic map in a neighbourhood of $0$. It follows from
assertion (i) in 
Corollary~\ref{maximum} that $f\circ v$ cannot have a local maximum
at $0$, unless $v$ is constant. In view of condition (h2) it
follows from assertion (ii) in Corollary~\ref{maximum} 
that for every $(s,t) \in \RR \times S^1$ for which
$u(s,t) \in \widehat{M} \setminus M$, the function $f \circ u$, which is 
well-defined in a neighbourhood of $(s,t)$, cannot have a local maximum
at $(s,t)$. But this contradicts the fact that $u(s,t)$ converges 
as $s \ra \infty$ to a periodic orbit which is entirely contained
in $M$. Hence $m = \lim_{s \to -\infty} u(s,t) \in M$. Now a similar 
reasoning as above, which uses again Corollary~\ref{maximum}, shows that
the whole image of $u$ lies in $M$. 
We have shown that $\cm = \cm_0$, 
and so Theorem~\ref{t:4finite} is proved. 
\proofend

Define the evaluation map $\mathrm{ev} \colon \cm \to M$ by
\[
\mathrm{ev}(u) \,:=\, \lim_{s \to -\infty} u(s,t) .
\]
Combining the techniques in \cite[Section 2.7]{Sa} and 
\cite[Appendix C.2]{F} one sees that the limit on the right-hand side
exists and 
that for generic $J_s$ the evaluation map $\mathrm{ev}$ is
transverse to every unstable manifold of the Morse function
$F \in C^\infty (M)$. 
Denote by $\Crit (F)$ the set of critical points of $F$ and by $\ind (c)$
the Morse index of $c \in \Crit (F)$.
Morse flow lines $\gg \colon \RR \ra M$ are solutions of the ordinary
differential equation
\begin{equation}  \label{e:4Morse}
\dot{\gg}(s) \,=\, - \nabla F(\gg(s)) 
\end{equation}
where the gradient is taken with respect to a generic
metric $g$ on $M$. 
For generators $c \in \Crit(F) \subset M$ of the Morse chain complex
and $x \in \cp_H$ of the Floer chain complex, let 
$\cm (c,x)$ be the moduli space of pairs $(\gg,u)$ such that 
$\gg \colon (-\infty ,0] \ra M$ solves \eqref{e:4Morse}, $u$ solves
\eqref{PSS}, and 
\[
\lim_{s \to -\infty} \gg (s) = c, 
 \quad\, 
\gg (0) = \mathrm{ev}(u) , 
\quad\, 
\lim_{s \to \infty} u(s,t) = x(t) .
\] 
If $\ind (c) = \mu(x)$, then $\cm (c,x)$
is a compact zero-dimensional manifold, see \cite{PSS}.
We can thus set
\[
n(c,x) \,:=\, \# \cm (c,x) \mod 2 .
\]
The Piunikhin--Salamon--Schwarz map $\phi \colon CM_*(M;F) \to
CF_*(M;H)$ is defined as the linear extension of
\[
\phi(c) = \sum_{\substack{c \in \cp_H\\ 
\ind (c)=\mu(x)}} n(c,x) \,x, \quad \, c \in \Crit(F) .
\]
By the usual gluing and compactness arguments
one proves that $\phi$ intertwines the boundary operators of
the Morse complex and the Floer complex and hence induces a homomorphism
\[
\Phi \colon HM_*(M) \,\ra\, HF_*(M) .
\]
To prove that $\Phi$ is an isomorphism we construct its inverse. 
We first define the Piunikhin--Salamon--Schwarz map 
$\psi \colon CF_*(M;H) \to CM_*(M;F).$ Let
\[
U \,:=\, \bigcup_{c \in \Crit (F)} W^s_F(c)
\]
be the union of the stable manifolds of $F$.
Since $F$ is admissible, the stable manifolds of the critical points 
of $F$ are entirely contained in the interior of $M$, i.e.,
\[
\overline{U} \subset M \setminus \pp M .
\]
For an open neighbourhood $V$ of $\overline{U}$ in $M \setminus \partial M$ 
choose a smooth family of admissible Hamiltonian functions $H_s$ for which
there exists $s_0>0$ such that
\[
H_s = H \;  \text{ if }\, s \leq -s_0 
\qquad \text{and} \qquad 
H_s|_V=0 \; \text{ if }\, s \geq s_0 .
\] 
Moreover, we assume that the Hamiltonian functions $H_s$ are the
restrictions of Hamiltonian functions 
$\widehat{H}_s \in \widehat{\ch}$, for which there exists 
$h \in C^\infty(\RR)$ independent of the $s$-variable which satisfies
$$h'(1)>0$$
such that
\[
\widehat{H}_s|_{S^1 \times \pp M \times [0,\infty)} \,=\, h \circ f .
\]
Choose a smooth family $J_s \in \cj (J^+,J^-)$
of admissible almost complex structures.
For $x \in \cp_H$ and $c \in \Crit (F)$ let
$\cm (x,c)$ be the moduli space of pairs $(u,\gg)$ 
such that $u$ solves \eqref{PSS},
$\gg \colon [0,\infty) \ra M$ solves \eqref{e:4Morse},
and
\begin{eqnarray*}
\lim_{s \ra -\infty} u(s,t) = x (t), 
 \quad\, 
 \lim_{s \ra \infty} u(s,t) = \gg (0) , 
 \quad\, 
\lim_{s \to \infty} \gg (s) = c .
\end{eqnarray*}
By our assumption on $\widehat{H}_s$ it follows from assertion (i) in
Corollary~\ref{maximum} that every solution of problem
\eqref{PSS2} is entirely contained in $M$ and hence solves
problem (\ref{PSS}). Hence we can show as above that
for generic choice of $J_s$ the moduli space $\cm (x,c)$
is a finite dimensional manifold of dimension
\[
\dim \cm (x,c)  \,=\, \mu (x) - \ind (c) .
\] 
In case that $\mu(x) = \ind (c)$, the moduli space is compact,
and we define
\[
n(x,c) \,:=\, \# \left\{ \cm(x,c) \right\} \mod 2.
\]
The Piunikhin--Salamon--Schwarz map $\psi \colon CF_*(M;H) \to CM_*(M;F)$
is defined as the linear extension of
\[
\psi(x) \, = \sum_{\substack{c \in \Crit (F)\\
\mu(x) = \ind (c)}} n(x,c) \,x , \quad\, x \in \cp_H .
\]
Again, $\psi$ intertwines the boundary
operators in the Floer complex and the Morse complex and hence
induces a homomorphism
\[
\Psi \colon HF_*(M) \,\ra\, HM_*(M) .
\]
One can prove that
\[
\Psi \circ \Phi = \id 
\qquad \text{and} \qquad 
\Phi \circ \Psi = \id,
\]
cf.\ \cite{PSS}, and so $\Phi$ and $\Psi$ are isomorphisms, called the
PSS isomorphisms.

\s
If $c_1 (\go)$ does not vanish on $\pi_2(M)$, we proceed in the same way
and obtain the PSS isomorphisms between the ungraded homologies
$HM (M)$ and $HF(M)$.
We refer to the appendix for a version of these isomorphisms preserving
a grading even if $c_1$ does not vanish on $\pi_2(M)$.

\subsection*{Products}
Proceeding as above and applying the maximum principle 
Corollary~\ref{maximum} factorwise we construct PSS isomorphisms also
for weakly exact compact split-convex symplectic manifolds. 

\section{The selector $c$}  \label{selector}

\ni
Let $(M, \go)$ be a weakly exact compact split-convex symplectic
manifold.
We do not assume that $c_1(\go)$ vanishes on $\pi_2(M)$
and shall work with ungraded chain complexes and homologies.
For a regular admissible Hamiltonian $H \in \ch_{\reg}$ 
and $a \in \RR$ let $CF^a(M;H)$ be the subvector space of
$CF(M;H)$ consisting of those formal sums
\[
\xi \,= \sum_{x \in \cp_H} \xi_x \,x, \quad\, \xi_x \in \mathbb{Z}_2 ,
\]
for which $\xi_x = 0$ if $\ca_H (x) > a$. 
In view of \eqref{est:action0}, the Floer boundary operator $\pp$
preserves $CF^a(M;H)$ and thus induces a boundary operator $\pp^a$
on the quotient $CF(M;H) / CF^a(M;H)$.
We denote the homology of the resulting complex by $HF^a(M;H)$. 
Since the projection $CF(M;H) \to CF(M;H)/CF^a(M;H)$ intertwines
$\pp$ and $\pp^a$, it induces a map 
$j^a \colon HF(M;H) \ra HF^a(M;H)$. 
Choose a generic admissible Morse function 
$F \in C^\infty(M)$ which attains its maximum in only one point, say $m$.
Let $[ \max ] \in HM(M)$ be the homology class represented by $m$.
Following \cite{Sch} we define
\begin{equation}  \label{def:cH}
c (H) \,:=\, \inf \left\{ a \in \RR \mid j^a \left( \Phi([ \max ]) \right) =
                                                0 \right\} 
\end{equation}
where $\Phi \colon HM (M) \to HF (M;H)$ is the PSS isomorphism.
Using the natural isomorphism $HF (M;H) \cong HF (M;K)$ for $H,K
\in \ch_{\reg}$ one can show that
\begin{equation}  \label{est:creg}
\left| c(H) - c(K) \right| \,\le\, \left\| H-K \right\|
\quad \text{ for all }\, H,K \in \ch_{\reg} ,
\end{equation}
see \cite[Section 2]{Sch}.
In particular, $c$ is $C^0$-continuous on $\ch_{\reg}$.
Let $\ch_c(M)$ be the set of $C^2$-smooth functions $S^1 \times M \ra
\RR$ whose support is contained in $S^1 \times \left( M \setminus \pp M
\right)$, and let $\ch_c^\infty (M)$ be the set of $C^\infty$-smooth
functions in $\ch_c(M)$.
Since $\ch_{\reg}$ is $C^\infty$-dense in $\ch$ 
and since $\ch_c^\infty (M)$ is $C^2$-dense in $\ch_c(M)$, 
we can first $C^\infty$-continuously extend $c$ to a map $\ch \to \RR$
and can then $C^2$-continuously extend its restriction to
$\ch_c^\infty(M)$ to a map $\ch_c(M) \ra \RR$ which we still denote by
$c$.
By \eqref{est:creg},
\begin{equation}  \label{cont:c:c}
\left| c(H) - c(K) \right| \,\le\, \left\| H-K \right\|
\quad \text{ for all }\, H,K \in \ch_c(M) .
\end{equation}
For $H \in \ch$ or $H \in \ch_c(M)$ we denote by $\cp_H$ the set of
contractible $1$-periodic orbits of $\gf^t_H$ and by $\Sigma_H$
the action spectrum   
\[
\Sigma_H \,=\, \left\{ \ca_H (x) \mid x \in \cp_H \right\} .
\]
The following property of $c$ is basic for everything to come.
\begin{proposition}  \label{p:cinspectrum}
For every $H \in \ch_c(M)$ it holds that $c(H) \in \Sigma_H$.
\end{proposition}

\proof
For $H \in \ch_{\reg}$ it follow from definition~\eqref{def:cH} that 
$c(H) \in \Sigma_H$.
For $H \in \ch_c(M)$ we choose a sequence $H_n$, $n \ge 1$, in
$\ch_{\reg}$ converging to $H$ in $C^2$ and choose 
$x_n \in \cp_{H_n}$ such that $c(H_n) = \ca_{H_n} (x_n)$.
Using that $M$ is compact we find a subsequence $n_j$, $j \ge 1$, 
such that $x_{n_j}(0) \ra x_0 \in M$ as $j \ra \infty$. 
Since the Hamiltonians $H_{n_j}$ converge to $H$ in $C^2$, it follows that 
$x(t) := \gf^t_H(x_0)$
belongs to $\cp_H$, and together with \eqref{cont:c:c},
\[
c(H) \,=\, \lim_{j \ra \infty} c \left( H_{n_j} \right)
     \,=\, \lim_{j \ra \infty} \ca_{H_{n_j}} \left( x_{n_j} \right)
     \,=\, \ca_H (x) . 
\]
Therefore, $c(H) \in \Sigma_H$.
\proofend

The set $\ch_c(M)$ forms a group with multiplication and inverse given by
\[
H_t \Diamond K_t = H_t + K_t \left( (\gf^t_{H_t})^- \right), 
\quad H_t^- = -H_t \circ \gf^t_{H_t}, 
\quad\, H_t, K_t \in \ch_c(M).
\]
It is shown in \cite{Sch} that $c$ satisfies the triangle inequality
\begin{equation}  \label{triangle}
c (H \Diamond K) \,\le\, c(H) + c(K), \quad\, H, K \in \ch_c(M) .
\end{equation}
The proof of \eqref{triangle} uses the product structure on Floer homology
given by the pair of pants product and a sharp energy estimate for the
pair of pants. 

\m
In the remainder of this section we give an upper bound for $c(H)$
and compute $c(H)$ for simple Hamiltonians.

\subsection{An upper bound for $c(H)$}  \label{action}
%
\begin{proposition}  \label{p:yaron} 
Let $(M, \go)$ be a weakly exact compact split-convex 
symplectic manifold, and let $H \in \ch_c(M)$.
Then
\begin{equation}  \label{e:cinf}
c(H) \,\le\, - \int_0^1 \inf_{x \in M} H_t (x) \, dt .
\end{equation}
In particular, $c(H)) \le \| H \|$.
\end{proposition}

\proof 
Since $c$ is $C^2$-continuous, it suffices to prove \eqref{e:cinf} for
$H \in \ch_{\reg}$.
Let $\widehat{H} \in \widehat{\ch}$ be such that $H = \widehat{H} |_{S^1
\times M}$.
We can choose the family $\widehat{H}_s \in C^\infty \big( S^1 \times
\widehat{M} \big)$ used in the construction of the PSS map
$\phi \colon CM (M;F) \ra CF (M;H)$ of the form
$\widehat{H}_s = \gb (s) \,\widehat{H}$ where $\gb \colon \RR \ra [0,1]$
is a smooth cut off function such that 
\begin{equation}  \label{es:beta}
\gb (s) =0,\,\, s \le 0;
\quad\,
\gb'(s) \ge 0,\,\, s \in \RR;
\quad\,
\gb (s) =1,\,\, s \ge 1 .
\end{equation}
In view of the construction of $\phi$ and the definition \eqref{def:cH} 
of $c(H)$ we find $x^+ \in \cp_H$ such that $\ca_H \left( x^+ \right) =
c(H)$ and a  solution $u \in C^\infty \left( \RR \times S^1, M \right)$
of the problem~\eqref{PSS}
such that 
$\lim_{s \ra \infty} u(s,t) = x^+(t)$.
Since the energy of $u$ is finite,
there exists $p \in M$ such that $\lim_{s \ra -\infty} u(s,t) =p$.
Using the Floer equation in \eqref{PSS} we compute
\begin{eqnarray*}
0 &\le& \int_0^1 \int_{-\infty}^\infty|\pp_s u|^2 \,ds \,dt\\
&=&-\int_0^1 \int_{-\infty}^\infty \left\langle \pp_s u,J_{s,t}(u) \left( \pp_t u-
X_{H_{s,t}}(u) \right) \right\rangle \,ds \,dt \\
&=&\int_{\RR \times S^1} u^*\go + \int_0^1 \int_{-\infty}^\infty
\go \left( X_{H_{s,t}}(u),\partial_s u \right) \,ds \,dt
\\
&=&\int_{\RR \times S^1}
 u^*\go + \int_0^1\int_{-\infty}^\infty d(H_{s,t}(u)) \pp_s u \,ds \,dt\\
&=&\int_{\RR \times S^1}
 u^*\go + \int_0^1\int_{-\infty}^\infty \frac{d}{ds} \Big( H_{s,t}(u) \Big)
\,ds \,dt\\
& &-\int_0^1 \int_{-\infty}^\infty \beta'(s) H_t(u) \,ds \,dt\\
&\le& \int_{\RR \times S^1}
 u^*\go + \int_0^1 H_t \left( x^+(t) \right) \,dt\\
& &- \left(\int_{-\infty}^\infty \beta'(s) \,ds \right) 
\left( \int_0^1 \inf_{x \in M} H_t(x) \,dt \right)\\
&=&-\ca_H \left( x^+ \right)-\int_0^1 \inf_{x \in M} H_t(x) \,dt .
\end{eqnarray*}
The proof of Proposition~\ref{p:yaron} is complete. 
\proofend

\subsection{A formula for $c(H)$.}
For a class of simple Hamiltonians the distinguished action value 
$c(H)$ can be explicitly computed.
The following theorem will be the main ingredient in the proof of the
energy-capacity inequality for the $\pi_1$-sensitive Hofer--Zehnder
capacity given in Section~\ref{ece}.

\begin{theorem}  \label{t:ece:cH}
Consider a weakly exact compact split-convex symplectic manifold $(M, \go)$,
and assume that $H \in \ch_c(M)$ has the following properties. 
\begin{itemize}
\item[(H1)] 
There exists $p \in \Int (M)$ such that
$H_t(p) =\min_{x \in M} H_t(x)$ for every $t \in [0,1]$.
\item[(H2)] 
The Hessian $\Hess (H)(p)$ of $H$ at $p$ with respect
to an $\go$-compatible Riemannian metric satisfies
\[
\left\| \Hess (H_t)(p) \right\| \,<\, 2 \pi
                        \quad\, \text{for all }\, t \in [0,1] .
\]
\item[(H3)] 
Every nonconstant periodic orbit of the flow $\gf_H^t$ has period 
greater than $1$. 
\end{itemize}
Then
\begin{equation}  \label{e:formula}
c(H) \,=\, -\int_0^1 H_t(p) \,dt.
\end{equation}
\end{theorem}

\proof 
It follows from assumptions (H1) and (H3) that the constant
orbit $p$ is a critical point of the action functional 
$\ca_{\lambda H}$ for every
$\lambda \in [0,1]$ and that for any other critical point $y$ of 
$\ca_{\lambda H}$,
\begin{equation}  \label{minact}
\ca_{\lambda H} (y) \,\le\, \ca_{\lambda H}(p)
\,=\, -\lambda \int_0^1 H_t(p) \,dt, \quad\, \lambda \in [0,1].
\end{equation}
We choose a sequence of regular admissible Hamiltonians 
$H_n \in \ch_{\reg}$ such that $H_n \ra H$ in $C^2$ 
and such that each $H_n$ satisfies (H1), (H2) and \eqref{minact} 
for the same point $p$. 
Since $c$ is $C^2$-continuous, it suffices to prove \eqref{e:formula}
for each $H_n$. We fix $n$ and from now on suppress $n$ in the notation.
We choose an admissible Morse function $F \in C^\infty(M)$ 
whose single maximum is attained at $p$,
and as in the previous paragraph we choose the family $H_s \in C^\infty
\left( S^1 \times M \right)$ of the form $H_s = \gb (s)\, H$ where $\gb
\colon \RR \ra [0,1]$ satisfies \eqref{es:beta}.
Let $c_p$ be the generator in $CM (M;F)$ represented by the maximum $p$
of $F$, and let $x_p$ be the generator of $CF(M;H)$ represented by $p$.
In view of the definition~\eqref{def:cH} of $c(H)$ and the construction
of the PSS map $\phi \colon CM (M;F) \ra CF (M;H)$, the formula
\eqref{e:formula} follows if we can show that
for generic choice of a smooth family $J_s$
of admissible almost complex structures which are independent of $s$
for $|s| \ge s_0$ large enough, the matrix coefficient
\[ 
n \left( c_p, x_p \right) \,=\, \# \cm \left( c_p , x_p \right) \mod 2
\] 
is odd.
Equivalently, we are left with showing
\begin{lemma}  \label{l:cob}
For generic choice of the smooth family $J_s$
of admissible almost complex structures independent of $s$
for $|s| \ge s_0$ large enough,
the number of solutions 
$u \in C^\infty( \RR \times S^1,M)$ of the problem 
\begin{equation}  \label{hom}
\left\{
  \begin{array}{rcl}
    \pp_s u + J_{s,t}(u) \left( \pp_t u-X_{H_{s,t}} (u) \right) 
                 &=& 0,  \\[0.2em]
    \displaystyle\lim_{s \to -\infty} u(s,t) &\in& W^u_F(p) , \\[0.4em]
    \displaystyle\lim_{s \to \infty}  u(s,t) &=& p,  \\[0.4em]
                c_1(u) &=& 0 ,                       
  \end{array}
 \right.
\end{equation}
is odd.
\end{lemma}

\proof
We choose a smooth family of
smooth families of admissible almost complex structures
$J^\lambda_s$, $s \in \RR$, $\lambda \in [0,1]$, such that 
$J^\lambda_s=J^{\lambda,\pm}$ is independent
of $s$ if $|s| \ge s_0$ is large enough, 
and consider for every $\lambda \in [0,1]$ the problem
\begin{equation}  \label{homhom}
\left\{
  \begin{array}{rcl}
\pp_s u+J^\gl_{s,t}(u) \left( \pp_t u-\gl X_{H_{s,t}}(u) \right) 
            &=& 0, \\[0.2em]
\displaystyle\lim_{s \to -\infty}u(s,t) &\in& W^u_F(p) , \\[0.2em]
\displaystyle\lim_{s \to \infty}u(s,t) &=& p, \\[0.2em]
     c_1(u) &=& 0 .
  \end{array}
 \right.
\end{equation}
Assumption (H2) guarantees that for each $\gl \in \;]0,1]$ the fixed point
$p$ of $\gf_{\gl H}^1$ is regular in the sense of \eqref{e:floer:det},
and hence for generic choice
of $J^\lambda_s$ the space $\cm_{\tot}$ of pairs $(u,\gl)$ solving
\eqref{homhom} for some $\gl \in [0,1]$ is 
a smooth $1$-dimensional manifold. 
The boundary $\pp \overline{\cm}_{\tot}$ of its compactification 
$\overline{\cm}_{\tot}$ contains an even number of elements,
\begin{equation}  \label{e:cob}
\# \pp \overline{\cm}_{\tot} \,=\, 0 \mod 2 .
\end{equation}
For generic choice of the family $J^\lambda_s$ transversality theory 
implies that $\pp \overline{\cm}_{\tot}$ consists of three types of
points, namely the solutions of \eqref{homhom} for $\gl =0$,
the solutions of \eqref{homhom} for $\gl =1$,
and broken trajectories.

\s
{\bf 1.} 
Since $\left[ \go \right]$ vanishes on $\pi_2(M)$,
the only solution of \eqref{homhom} for $\lambda=0$ 
is the constant map $u \equiv p$. 

\s
{\bf 2.} 
The solutions of \eqref{homhom} for $\gl =1$ are the solutions of
\eqref{hom} which we want to count.

\s
{\bf 3.} 
Solutions of \eqref{homhom} are in bijection with solutions 
consisting of half a Morse flow line followed by a Floer
disc.
For generic choice of the family $J_s^\gl$, these solutions break off
only once, either along the Morse flow line or along the Floer disc.
More precisely, for generic choice of $J_s^\gl$, 
there are finitely many values
$0 < \gl_1 < \ldots < \gl_n <1$ for which there are broken trajectories
consisting either of pairs 
$u_1 \in C^\infty(\RR,M)$,  $u_2 \in C^\infty(\RR \times S^1,M)$ which
satisfy, for some $i \in \left\{ 1, \dots, n \right\}$, 
\begin{equation}  \label{left}
\left\{
  \begin{array}{rcl}
\pp_s u_1 &=& -\nabla F(u_1), \\[0.2em]
\pp_s u_2+J^{\gl_i}_{s,t} \left( \pp_t u_2-\gl_iX_{H_{s,t}} (u_2)
                                        \right)  &=& 0, \\[0.2em]
\displaystyle\lim_{s \to -\infty}u_1(s,t) &=& p, \\[0.4em]
\displaystyle\lim_{s \to \infty}u_1(s,t)  &\in& \Crit \left(F\right), \\[0.4em]
\ind_F \Big( \displaystyle\lim_{s \to \infty}u_1(s,t) \Big) &=& 2n-1, \\[0.4em]
\displaystyle\lim_{s \to -\infty} u_2(s,t) &\in& W^u_F \left(
\displaystyle\lim_{s \to \infty} u_1(s) \right), \\[0.4em]
\displaystyle\lim_{s \to \infty} u_2(s) &=& p, \\[0.4em]
c_1(u_2) &=& 0 ,
 \end{array}
 \right.
\end{equation}
or pairs
$u_1, u_2 \in C^\infty (\RR \times S^1,M)$
which satisfy, for some $i \in \left\{ 1,\ldots,n \right\}$,
\begin{equation}  \label{right}
\left\{
  \begin{array}{rcl}
\pp_s u_1+J^{\gl_i}_{s,t} \left( \pp_t u_1-\gl_iX_{H_{s,t}}(u_1) \right) &=&
0, \\[0.2em] 
\pp_s u_2+J^{\gl_i,+}_t \left( \pp_t u_2-\gl_i X_{H_t}(u_2)
\right) &=&0, \\[0.2em]
\displaystyle\lim_{s \to -\infty} u_1(s,t) &\in& W^u_F(p), \\[0.4em]
\displaystyle\lim_{s \to \infty} u_1(s,t) \,=\, \lim_{s \to -\infty}u_2(s,t)
&\in& \Crit \left( \ca_{\gl_i H} \right) \setminus \{p\} ,\\[0.4em]
\displaystyle\lim_{s \to \infty} u_2(s,t) &=& p , \\[0.4em]
           c_1(u_1\#u_2) &=& 0 ,
  \end{array}
 \right.
\end{equation}
where the sphere $u_1 \# u_2$ is the connected sum of the oriented discs
$u_1$ and $u_2$. 
Since $p$ is the only maximum of $F$, 
for each critical point of $F$ of index $2n-1$ there is an even number
of Morse flow lines $u_1$ emanating from that point and ending in $p$. 
This shows that there is an even number of solutions of \eqref{left}.
Moreover, it follows from formula~\eqref{est:action0} and from
assumption~\eqref{minact} that solutions $u_2$ of problem 
\eqref{right} have nonpositive energy and hence cannot exist.
We conclude that there is an even number of broken trajectories.

\s
In view of \eqref{e:cob} and 1.\ and 3.\ we conclude that for generic
choice of $J_s$ the number of solutions of \eqref{hom} is odd.
This proves Lemma~\ref{l:cob}, and so Theorem~\ref{t:ece:cH} is also proved.
\proofend

\section{The action spectrum}

\ni
Recall that the action spectrum $\Sigma_H$ of $H \in \ch_c(M)$ is the set
\[
\Sigma_H \,=\, \left\{ \ca_H (x) \mid x \in \cp_H \right\} .
\]
For a {\it closed}\, symplectic manifold, the dependence of the action
spectrum on the Hamiltonian $H$ is a subtle problem, see \cite{Oh,Sch}.
As we shall see in this section, for an open (i.e., not closed)
weakly exact symplectic manifold, $\Sigma_H = \Sigma_K$ whenever 
$H,K \in \ch_c(M)$ generate the same Hamiltonian diffeomorphisms 
$\gf_H = \gf_K$.

Let $(M, \go)$ be an open weakly exact symplectic manifold, and
let $G \in \ch_c(M)$ be such that $\gf_G = \id$. 
To $q \in M$ we associate the loop 
\[
x_q(t) := \gf_G^t(q), \quad t \in [0,1] . 
\]
If $q \in M \setminus \supp \gf_G$, then $x_q$ is the constant loop.
This and the continuity of the map $q \mapsto x_q$ from $M$ to the free
loop space of $M$ show that $x_q \in \cp_G$ for all $q \in M$.
We define the function $I_G \colon M \ra \RR$ by
\[
I_G(q) \,\equiv\, \ca_G (x_q) \,=\, 
- \int_{D^2} \bar{x}_q^* \go - \int_0^1 G (t,x_q(t) )\, dt
\]
where $\bar{x}_q$ is a smooth extension of $x_q$ to the unit disc $D^2$.
\begin{proposition}  \label{monodromy}
The function $I_G$ vanishes identically.
\end{proposition}

\proof 
 If $q \in M \setminus \supp \gf_G$, then $I_G (q) =0$.
It remains to show that $I_G$ is constant.
To this end we choose a path $r \mapsto q(r)$ and compute
\begin{eqnarray*}
\frac{d}{dr} I_G (q(r)) &=& - \int_0^1 \go \left( d \gf^t_G (q) q'(r),
                             X_{G_t} \left( \gf^t_G (q) \right) \right) dt \\
                        & & - \int_0^1 dG_t (\gf^t_G(q)) \left( d
                                           \gf^t_G(q) q'(r) \right) dt
                                           \;=\; 0, \\
\end{eqnarray*}
as desired.
\proofend

Consider $H,K \in \ch_c(M)$ such that $\gf_H = \gf_K$.
We choose a smooth function $\ga \colon [0,1] \ra [0,1]$ such that
\begin{equation}  \label{e:alpha}
\ga (t) \,=\, 
 \left\{
  \begin{array}{ll}
   0, & t \le 1/6 , \\ [0.2em]
   1, & t \ge 1/3 .
  \end{array}
 \right.
\end{equation}
The Hamiltonian $G \in \ch_c(M)$ defined by
\begin{equation}  \label{e:G}
G (t,x) \,=\, 
 \left\{
  \begin{array}{ll}
    \ga'(t)  \,H(\ga(t),x),   & 0 \le t \le 1/2 , \\ [0.4em]
   -\ga'(1-t)\,K(\ga(1-t),x), & 1/2 \le t \le 1 ,
  \end{array}
 \right.
\end{equation}
then generates the loop
\begin{equation*}  
\gf_G^t \,=\, 
 \left\{
  \begin{array}{ll}
   \gf_H^{\ga(t)},   & 0 \le t \le 1/2 , \\ [0.4em]
   \gf_K^{\ga(1-t)}, & 1/2 \le t \le 1 ,
  \end{array}
 \right.
\end{equation*}
in $\Ham^c(M, \go)$.
Since all loops $x_q(t) = \gf_G^t(q)$, $q \in M$, $t \in [0,1]$,
are contractible, the sets $\cp_H$ and $\cp_K$ can be canonically
identified, and the set
\[
\Fix (\gf_H) \,=\, \left\{ x(0) \mid x \in \cp_H \right\}
\]
of ``contractible fixed points'' of $\gf_H$ does not depend on $H$.
The action of a fixed point $x \in \Fix (\gf_H)$ is defined as the
action of the loop $\gf_H^t(x)$,
\[
\ca_H (x) \,:=\, \ca_H \left( \gf_H^t(x) \right) .
\]
\begin{corollary}  \label{c:indep}
Assume that $H,K \in \ch_c(M)$ are such that $\gf_H = \gf_K$.
Then
$\ca_H (x) = \ca_K(x)$ for all $x \in \Fix (\gf_H)$.
In particular, $\Sigma_H = \Sigma_K$.
\end{corollary}

\proof
Define $G \in \ch_c(M)$ as in \eqref{e:G}. Then
\[
\ca_G \left( \gf_G^t (x) \right) \,=\, 
\ca_H \left( \gf_H^t (x) \right) - \ca_K \left( \gf_K^t (x) \right)
\]
for all $x \in \Fix (\gf_H) = \Fix (\gf_K)$, and so
Corollary~\ref{c:indep} follows from Proposition~\ref{monodromy}.
\proofend

Recall that the inverse of $H \in \ch_c(M)$ is defined as
\[
H^-_t (x) := - H_t \left( \gf^t_H(x) \right) .
\]
To $x \in \cp_H$ we associate the loop $x^-$ defined as
\[
x^-(t) := \gf^t_{H^-}(x(0)) .
\]
\begin{corollary}  \label{c:inv}
If $x \in \cp_H$, then $x^- \in \cp_{H^-}$ and 
$\ca_{H^-} (x^-) = -\ca_H(x)$.
In particular, $\Sigma_{H^-} = - \Sigma_H$.
\end{corollary}

\proof 
Choose $\ga \colon [0,1] \ra [0,1]$ as in \eqref{e:alpha} 
and define $G \in \ch_c(M)$ by
\begin{equation*} 
G (t,x) \,=\, 
 \left\{
  \begin{array}{ll}
    \ga'(t)  \,H(\ga(t),x),         & 0 \le t \le 1/2 , \\ [0.4em]
    \ga'(t-1/2)\,H^-(\ga(t-1/2),x), & 1/2 \le t \le 1 .
  \end{array}
 \right.
\end{equation*}
Then $\gf_G=\id$. For $x \in \cp_H$ the loop $x^-$ therefore belongs to
$\cp_{H^-}$. Moreover, 
\[
I_G(x(0)) \,=\, \ca_H(x) + \ca_{H^-} (x^-) ,
\]
and so Proposition~\ref{monodromy} yields $\ca_{H^-} (x^-) = -\ca_H(x)$.
Since the map $x \mapsto x^-$ is a bijection between $\cp_H$ and
$\cp_{H^-}$, we conclude $\Sigma_{H^-} = - \Sigma_H$.
\proofend

\section{The Schwarz metric}  \label{metric}

\ni
We consider a weakly exact compact split-convex
symplectic manifold $(M, \go)$.
For $H \in \ch_c(M)$ let $c(H) \in \Sigma_H$ be the the
distinguished critical value of $\ca_H$ defined in Section~\ref{selector}.
\begin{proposition}  \label{p:metric:indep}
Assume that $H,K \in \ch_c(M)$ satisfy $\varphi_H = \varphi_K$. Then
\[
c(H) = c(K).
\]
\end{proposition} 
\proof 
Since $\gf_{H \Diamond K^-} = \gf_0 = \id$, 
Corollary~\ref{c:indep} shows
that $\Sigma_{H \Diamond K^-} = \Sigma_0 = \{ 0 \}$.
This and Proposition~\ref{p:cinspectrum} yield $c \left( H \Diamond K^-
\right) =0$. Together with the triangle inequality~\eqref{triangle}
we conclude
\[
c \left( H \right) \,=\, c \left( H \Diamond K^- \Diamond K \right) \,\le\, 
c \left( H \Diamond K^- \right) + c \left( K \right) \,=\, 
c \left( K \right) .
\]
Interchanging the roles of $H$ and $K$ we
obtain $c(K) \le c(H)$.
Proposition~\ref{p:metric:indep} follows. 
\proofend

In view of Proposition~\ref{p:metric:indep} we can define $c \colon
\Ham_c(M) \ra \RR$ by  
\[
c(\gf) = c(H) \quad \text{if }\, \gf = \gf_H .
\]
We define the {\it Schwarz norm} $\gg \colon \Ham_c(M) \ra \RR$ by
\begin{equation}  \label{def:gamma}
\gg \left( \gf \right) \,=\, c \left( \gf \right) + c \left( \gf^{-1}
\right) .
\end{equation}
We shall often write $\gg (H)$ instead of $\gg \left( \gf_H \right)$.
By Proposition~\ref{p:cinspectrum} and Corollary~\ref{c:inv},
$c(H) \in \Sigma_H$ and $- c \left( H^- \right) \in \Sigma_H$, and so
$\gg \left( \gf_H \right) = \gg \left( H \right) = c \left( H \right) + c
\left( H^- \right)$ is the difference of two distinguished actions of $\gf_H$.
\begin{proposition}  \label{p:C2small}
For every $C^2$-small time-independent $H \in \ch_c(M)$ we have
$
\gg (H) = \left\| H \right\| .
$
\end{proposition}

\proof
According to Theorem~\ref{t:ece:cH} we have $c(H) = - \min H$ and $c
\left( H^- \right) = c \left( -H \right) = \max H$, and so $\gg \left( H
\right) = c \left( H \right) + c \left( H^- \right) = \left\| H
\right\|$.
\proofend

We recall that $\Sympcc (M)$ denotes the group of symplectomorphisms 
of $(M, \go)$ whose support is contained in $M \setminus \pp M$.
The following theorem justifies that $\gg$ is called a norm.
\begin{theorem}  \label{t:norm}
The Schwarz norm $\gg$ on $\Ham_c(M)$ has the following properties.
\begin{itemize}
\item[(S1)]
$\gg (\idd) =0$ and $\gg (\gf) >0$ if $\gf \neq \idd$;
\item[(S2)]
$\gg (\gf \psi) \le \gg (\gf) + \gg (\psi)$;
\item[(S3)]
$\gg ( \gt \gf \gt^{-1}) = \gg (\gf)$ for all $\gt \in \Sympcc (M)$;
\item[(S4)]
$\gg (\gf) = \gg \left( \gf^{-1} \right)$;
\item[(S5)]
$\gg (\gf) \le d_H \left( \gf, \idd \right)$.
\end{itemize}
\end{theorem}

\proof
The triangle inequality (S2) follows from the triangle
inequality~\eqref{triangle} for $c$.
For $\gf_H \in \Ham_c (M)$ and $\gt \in \Sympcc (M)$ we have 
\[
\gt \circ \gf_H^t \circ \gt^{-1} \,=\, \gf_{H_\gt}^t \quad\, 
\text{for all }\, t
\]
where $H_{\gt} (t,x) = H \left( t, \gt^{-1}(x) \right)$.
This and the invariance of the Floer equation
imply the invariance property (S3).
The symmetry property (S4) follows from definition~\eqref{def:gamma}.
In order to prove the estimate (S5) we need to show that
$c \left( H \right) + c \left( H^- \right) \le \left\| H \right\|$ for
all $H \in \ch_c(M)$.
In view of the continuity of $c$, it suffices to show this for $H \in
\ch_{\reg}$.
According to Proposition~\ref{p:yaron} we have
\begin{equation}  \label{e:cHinf}
c \left( H \right) \,\le\, - \int_0^1 \inf_{x \in M} H_t(x) \,dt ,
\end{equation}
and combining Proposition~\ref{p:yaron} with
\[
\inf_{x \in M} H_t^- (x) \,=\,
\inf_{x \in M} \left( - H_t \left( \gf_H^t (x) \right) \right) \,=\,
\inf_{x \in M} \left( -H_t(x) \right) \,=\,
- \sup_{x \in M} H_t(x) 
\]
we find
\begin{equation}  \label{e:cHsup}
c \left( H^- \right) \,\le\, - \int_0^1 \inf_{x \in M} H_t^-(x) \,dt 
\,=\, \int_0^1 \sup_{x \in M} H_t(x)\, dt .
\end{equation}
Adding \eqref{e:cHinf} and \eqref{e:cHsup} we obtain $c \left( H \right)
+ c \left( H^- \right) \le \left\| H \right\|$, as desired.
\proofend

We are left with proving (S1).
If $\gf = \gf_0 = \id$, then 
$c \left( \gf \right) = c \left( \gf^{-1} \right) = c \left( 0 \right)
=0$ 
and so $\gg \left( \gf \right) = 0$.
In order to verify that $\gg$ is non-degenerate, we shall need the
following proposition, which will be crucial for most of our
applications.
\begin{proposition}  \label{p:phi2psi}
Assume that $\gf_H, \psi \in \Ham_c(M)$ are such that $\psi$ displaces
$\supp \gf_H$.
Then $\gg \left( \gf_H^n \right) \le 2\, \gg \left( \psi \right)$
for all $n \in \NN$.
\end{proposition}

\proof
We closely follow the proof of Proposition~5.1 in \cite{Sch}.
Since $\supp \gf_H^n = \supp \gf_H$ for all $n \in \NN$, it is enough to
prove the claim for $n=1$.
Assume that $\psi = \gf_K$. After reparametrizing in $t$ we can assume
that $H_t=0$ for $t \in \left[ 0, 1/2 \right]$ and $K_t=0$ for 
$t \in \left[ 1/2, 1 \right]$.
With this choice of $H$ and $K$ and since $\gf_K$
displaces $\supp \gf_{\eps H} = \supp \gf_H$ for each $\eps \in [0,1]$
it is clear that 
\[
\Fix \left( \gf_{\eps H \Diamond K} \right) \,=\, 
\Fix \left( \gf_K \right) \,\subset\,
M \setminus \supp \gf_{\eps H} ,
\]
and so $\cp_{\eps H \Diamond K} = \cp_K$ and 
$\Sigma_{\eps H \Diamond K} = \Sigma_K$ for each $\eps \in [0,1]$.
The set $\Sigma_K = \Sigma_{\eps H \Diamond K}$ is nowhere dense, see
\cite[Proposition 3.7]{Sch}.
This and the continuity of $c$ imply that the map
\[
[0,1] \ra \Sigma_K, \quad\, 
\eps \mapsto c \left( \eps H \Diamond K \right) ,
\]
is constant.
In particular, $c \left( H \Diamond K \right) = c \left( K \right)$.
Since $\gf_K$ displaces $\supp \gf_H$, its inverse $\gf_{K^-}$ displaces
$\supp \gf_{H^-} = \supp \gf_H$.
An argument analogous to the above then yields 
$c \left( \left(H \Diamond K \right)^- \right) = c \left( K^- \Diamond
H^- \right) = c \left( K^- \right)$.
Summarizing we find $\gg \left( H \Diamond K \right) = \gg \left( K
\right)$.
Together with (S2) and (S4) we can thus conclude
\[
\gg \left( H \right) \,=\, 
\gg \left( H \Diamond K \Diamond K^- \right) \,\le\, 
\gg \left( H \Diamond K \right) +\gg \left( K^- \right) \,=\, 
2 \,\gg \left( K \right) ,
\]
as desired.
\proofend

Assume now that $\gf \neq \id$. We then find a non-empty open subset 
$U \subset M$ such that $\gf$ displaces $U$.
According to Proposition~\ref{p:C2small} we can choose $H \in \ch_c(M)$
such that $\gg (H) >0$.
Applying Proposition~\ref{p:phi2psi} with $\psi = \gf$ we get $0< \gg (H)
\le 2 \gg (\gf)$.
The proof of Theorem~\ref{t:norm} is complete.
\proofend

\begin{corollary}  \label{c:spectrumnontrivial}
If $\gf_H \in \Ham_c(M) \setminus \{ \idd \}$, then the spectrum
$\Sigma_H$ contains not only $0$.
\end{corollary}

\proof
Recall that $\gg (H)$ is the difference of two elements of $\Sigma_H$.
The corollary thus follows from (S1) of Theorem~\ref{t:norm}.
\proofend

The Schwarz metric on $\Ham_c(M)$ is defined as
\[
d_S \left( \gf,\psi \right) \,:=\, 
\gg \left( \gf \circ \psi^{-1} \right) \quad\, \gf, \psi \in \Ham_c (M).
\]
Theorem~\ref{t:norm} says that $d_S$ is a biinvariant metric on
$\Ham_c(M)$ such that
\[
d_S \left( \gf, \psi \right) \,\le\, d_H \left( \gf, \psi \right) \quad\,
\text{for all }\m \gf, \psi \in \Ham_c(M) .
\]

\section{An energy-capacity inequality}  \label{ece}

\ni
In this section we shall compare the $\pi_1$-sensitive Hofer--Zehnder
capacity $c_{\HZ}^\circ (A)$ of a subset $A \subset (M, \go)$ with the
Schwarz-diameter of $\Ham_c \left( \Int A, \go \right)$. This will lead
to an energy-capacity inequality for $c_{\HZ}^\circ$,
which will be a crucial tool in the proofs of Theorems~4.A and 4.B\:(ii).

\subsection{The $\pi_1$-sensitive Hofer--Zehnder capacity.}

Let $(M, \go)$ be an arbitrary symplectic manifold.
Given a subset $A \subset M$ we consider the function space
\[
\ch_c(A) \,=\, \left\{ H \in C^\infty_c \left( \Int A \right) \mid
        H \ge 0,\, H |_U = \max H \text{ for some open } U
        \subset A \right\} .
\]
We say that $H \in \ch_c(A)$ is {\it $\HZ$-admissible}\, if the flow
$\gf_H^t$ has no non-constant $T$-periodic orbit with period $T \le 1$,
and we say that $H \in \ch_c(A)$ is {\it $\HZ^\circ$-admissible}\, 
if the flow $\gf_H^t$ has no non-constant $T$-periodic orbit 
with period $T \le 1$ which is contractible in $M$.
Set
\begin{eqnarray*}
\ch_{\HZ} (A,M,\go) &=& \left\{ H \in \ch_c(A) \mid H \text{ is
HZ-admissible} \right\}, \\[0.2em]
\ch_{\HZ}^\circ (A,M,\go) &=& \left\{ H \in \ch_c(A) \mid H \text{ is
$\HZ^\circ$-admissible} \right\}.
\end{eqnarray*}
The Hofer--Zehnder capacity and the $\pi_1$-sensitive Hofer--Zehnder
capacity of $A \subset (M, \go)$ are defined as
\begin{eqnarray*}
c_{\HZ} (A,M,\go) &=& \sup \left\{ \left\| H \right\| \mid H \in \ch_{\HZ}
(A,M,\go) \right\}, \\[0.2em]
c_{\HZ}^\circ (A,M,\go) &=& \sup \left\{ \left\| H \right\| \mid H \in
\ch_{\HZ}^\circ (A,M,\go) \right\} .
\end{eqnarray*}
From now on we suppress $\go$ from the notation. Of course, $c_{\HZ}
(A,M) \le c_{\HZ}^\circ (A,M)$.
Example~\ref{ex:cap} below shows that this inequality can be strict.
It also shows that in contrast to $c_{\HZ}$, the $\pi_1$-sensitive
Hofer--Zehnder capacity $c_{\HZ}^\circ$ is not an intrinsic symplectic
capacity as defined in \cite{HZ};
it is, however, a relative symplectic capacity and in particular
satisfies the relative monotonicity axiom
\begin{equation}  \label{e:relmon}
c_{\HZ}^\circ (A,M) \,\le\, c_{\HZ}^\circ (B,M) \quad
\text{ whenever } A \subset B \subset M .
\end{equation}

\begin{example}  \label{ex:cap}
{\rm
Consider the annulus $A = \left\{ z \in \RR^2 \mid 0<\left| z \right| <1
\right\}$ in $\left( \RR^2, \go_0 \right)$.
Then
$c_{\HZ} (A,A) = c_{\HZ}^\circ (A,\RR^2) = \pi$ and
$c_{\HZ}^\circ (A,A) = \infty$.
}
\end{example}

\begin{corollary}  \label{c:cgM}
For any subset $A$ of a weakly exact compact split-convex symplectic manifold 
$(M, \go)$, 
\[
c_{\HZ}^\circ (A,M) \,=\, \sup \left\{ \gg_M \left( \gf_H \right) \mid
H \in \ch_{\HZ}^\circ (A,M) \right\} .
\]
\end{corollary}  

\proof
Fix $H \in \ch_{\HZ}^\circ (A,M)$. Then both $H$ and $H^- = -H$ meet the
assumptions of Theorem~\ref{t:ece:cH}, and so
$c(H) =0$ and $c \left( H^- \right) = \left\| H \right\|$.
Therefore, $\gg_M(H) = c(H)+c \left( H^- \right) = \left\| H \right\|$.
\proofend

\subsection{An energy-capacity inequality for $c_{\HZ}^\circ$.}

Following \cite{Sch} we define for any subset $A$ of a weakly exact
compact split-convex symplectic manifold $(M,\go)$ the relative
capacity $c_\gg (A,M) = c_{\gg} (A,M,\go) \in [0,\infty]$ as
\[
c_\gg (A,M) \,=\, \sup \left\{ \gg_M (\gf) \mid \gf \in \Ham_c (M,
\go), \, \supp \gf \subset S^1 \times A \right\} . 
\]
Notice that $c_\gg (A,M)$ is the diameter of $\Ham_c \left( \Int A
\right\}$.
We recall that the displacement energy $e(A,M) = e (A,M,\go)$ is defined
as
\[
e(A,M) \,=\, \inf \left\{ d_H (\gf, \id ) \mid \gf \in \Ham_c (M, \go),
\, \gf(A) \cap A = \emptyset \right\} .
\]
Corollary~\ref{c:cgM}, Proposition~\ref{p:phi2psi} and (S5) of
Theorem~\ref{t:norm} yield 
\begin{corollary}  \label{c:ece}
For any subset $A$ of a weakly exact compact split-convex
symplectic manifold $(M,\go)$, 
\[
c_{\HZ} (A,M) \,\le\, c_{\HZ}^\circ (A,M) \,\le\, c_\gg (A,M) \,\le\,  
2 \,e (A,M) .
\]
\end{corollary}  

\m
This concludes the construction of our tools. In the next five sections
we shall use them to study Hamiltonian diffeomorphisms on weakly exact
symplectic manifolds which away 
from a compact subset look like a product of convex symplectic manifolds.
To be precise,
we recall from Definition~\ref{def:splitconvex} that 
a compact symplectic manifold $(M,\go)$ is split-convex
if there exist compact convex symplectic manifolds 
$(M_j, \go_j)$, $j = 1, \dots, k$, and a compact subset 
$K \subset M \setminus \pp M$ such that $M = M_1 \times \dots \times
M_k$ and 
\[
\left( M \setminus K, \go \right) \,=\, \left( \left( M_1 \times
\dots \times M_k \right) \setminus K, \go_1 \oplus \dots \oplus \go_k
\right) .
\]

We say that a non-compact symplectic manifold $(M, \go)$ is 
{\it split-convex}\, if 
there exists an increasing sequence of compact split-convex submanifolds 
(with corners) $M_i \subset M$ exhausting $M$, that is,
\[
M_1 \subset M_2 \subset \dots \subset M_i \subset \dots \subset M
\quad \text{ and } \quad
\bigcup_i M_i = M . 
\]  

\section{Existence of a closed orbit with non-zero action}  \label{growth}

\ni
The following result is a generalization of Theorem~1.

\begin{theorem}  \label{t:action:products}
Assume that $(M, \go)$ is a weakly exact split-convex symplectic
manifold. 
Then for every Hamiltonian function $H \in \ch_c(M)$ generating 
$\gf_H \in \Ham_c(M,\go) \setminus \left\{ \idd \right\}$ there
exists $x \in \cp_H$ such that $\ca_H(x) \neq 0$.
\end{theorem}

\proof
Assume that $\bigcup_{i \ge 1} M_i$ is an exhaustion of $M$ by compact
split-convex submanifolds.
Given $H \in \ch_c(M)$ generating $\gf_H \neq \id$
we choose $i$ so large that $\supp \gf_H \subset M_i$.
Since $\gf_H \in \Ham_c \left(M_i \right) \setminus \left\{ \idd
\right\}$, Corollary~\ref{c:spectrumnontrivial} guarantees the existence
of $x \in \cp_H$ with $\ca_H (x) \neq 0$,
and so Theorem~\ref{t:action:products} follows.
\proofend

\section{Infinitely many periodic points of Hamiltonian diffeomorphisms}

\ni
We first consider a weakly exact compact split-convex 
symplectic manifold $(M,\go)$,
and we let $\gg$ be the Schwarz norm on $\Ham_c(M, \go)$ constructed
in Section~\ref{metric}.

\begin{theorem}  \label{c:bound}
Assume that $\gf_H \in \Ham_c (M, \go) \setminus \left\{ \idd \right\}$ 
is such that
\[
\text{
$\gg (\gf_H^n) \le C$\, for all $n \in \NN$ and some $C<\infty$.
}
\]
Then $\gf_H$ has infinitely many nontrivial geometrically distinct
periodic points corresponding to contractible periodic orbits.  
\end{theorem}  

\proof
We closely follow \cite{Sch}.

\s
\ni
{\bf Case 1.} {\it $\gf_H^n = \idd$ for some $n \in \NN$.}
Then every $x \in M$ is a periodic point of $\gf_H$, and since the
support of $\gf_H$ is not all of $M$ and since $M$ is connected, every
$x \in M$ is a periodic point of $\gf_H$ corresponding to a contractible
periodic orbit.
Since $\gf_H \neq \idd$, infinitely many among these periodic points are
non-trivial.

\s
\ni
{\bf Case 2.} {\it $\gf_H^n \neq \idd$ for all $n \in \NN$.}
According to Corollary~\ref{c:spectrumnontrivial}, $\gf_H$ has at least
$1$ nontrivial periodic point corresponding to a contractible
periodic orbit.
Arguing by contradiction, we assume that $\gf_H$ has only finitely many
nontrivial geometrically distinct periodic points corresponding to
contractible periodic orbits, say $x_1, \dots, x_N$. 
The period of $x_i$ is defined as the minimal
$k_i \in \NN$ such that $\gf^{k_i}_H(x_i) = x_i$.
Set $k = k_1 k_2 \cdots k_N$ and $G(t,x) =  k H (kt, x)$.
Then $\gf_G = \gf_H^k$, and $x_1, \dots, x_N$ are the nontrivial
periodic points of $\gf_G$ corresponding to contractible periodic
orbits.
There period is $1$.
By assumption,
\begin{equation}  \label{est:C}
\text{
$\gg \left( \gf_G^n \right) \,=\, \gg \left( \gf_H^{nk} \right) \,\le\,
C$\, for all $n \in \NN$. 
}
\end{equation}
The spectrum $\Sigma_G$ consists of $0$ (coming from trivial periodic
points) and $\ca_G(x_i)$, $i = 1, \dots, N$.
Set $G^{(n)} (t,x) = n G(nt,x)$. 
Since $\gf_G$ has no other nontrivial periodic points 
corresponding to contractible periodic orbits than 
$x_1, \dots, x_N$,
\begin{equation}  \label{id:specn}
\Sigma_{G^{(n)}} \,=\, n \Sigma_G \,=\,
\left\{ 0, n \ca_G(x_1), \dots, n \ca_G(x_N) \right\} .
\end{equation}
By assumption, $\gf_G^n = \gf_H^{nk} \neq \id$ for all $n$, and so
\[
\gg \left( \gf_G^n \right) \,=\, \gg \left( G^{(n)} \right) \,=\, 
c \left( G^{(n)} \right) + c \,\big( \left( G^{(n)} \right)^- \big)
\,>\, 0 
\quad\,\text{for all }\, n \in \NN . 
\]
Recall now that 
$c \left( G^{(n)} \right) + c \, \big( \left( G^{(n)} \right)^- \big)$ 
is the difference of two action values in $\Sigma_{G^{(n)}}$.
We thus infer from \eqref{id:specn} that 
$\gg \left( \gf_G^n \right) \ra \infty$ as $n \ra \infty$, contradicting
\eqref{est:C}.
\proofend

Theorem~2 is a special case of 
\begin{corollary}  \label{c:displaceable}
Assume that $(M, \go)$ is a weakly exact split-convex symplectic manifold.
If the support of 
\text{$\gf_H \in \Ham_c (M, \go) \setminus \{ \idd \}$}
is displaceable,
then $\gf_H$ has infinitely many nontrivial
geometrically distinct periodic points corresponding to contractible 
periodic orbits.  
\end{corollary}

\proof
Choose $\psi \in \Ham_c(M, \go)$ which displaces $\supp \gf_H$, and
choose $i$ so large that $\supp \psi \subset M_i$.
According to Proposition~\ref{p:phi2psi},
$\gg_{M_i} \left( \gf_H^n \right) \le 2\, \gg_{M_i} \left( \psi \right)$
for all $n \in \NN$,
and so the corollary follows from Theorem~\ref{c:bound}.
\proofend

\ni
{\it Proof of Corollary 2:}
Consider a subcritical Stein manifold $(V, J, f)$ and 
$\gf_H \in \Ham_c \left( V, \go_f \right) \setminus \{ \idd \}$.
Since $f$ is proper, we find a regular value $R$ such that 
$S = \supp \gf_H$ 
is contained in $V_R = \left\{ x \in V \mid f(x) \le R \right\}$.
After composing $f$ with an appropriate smooth function
$h \colon \RR \ra \RR$ such that $h(r) = r$ for $r \le R$
we obtain a subcritical Stein manifold $(V,J, h \circ f)$
such that the gradient vector field $X_{h \circ f}$ of $h \circ f$ with
respect to the Riemannian metric $g_{h\circ f}$ is complete, 
see \cite[Lemma 3.1]{BC}.
Since $S \subset V_R$ and $\go_f |_{V_R} = \go_{h \circ f} |_{V_R}$, we
have $\gf_H \in \Ham_c \left( V, \go_{h\circ f} \right) \setminus \{
\idd \}$. 
Let $\Crit_R (h \circ f)$ be the set of critical points of $h
\circ f$ in $V_R$, and consider the union
\[
\Delta_R \,=\, 
\bigcup_{x \in \Crit_R (h \circ f)} W_x^s \left( X_{h \circ
f} \right) 
\]
of those stable submanifolds of $X_{h\circ f}$ which are contained in $V_R$. 
Applying the proof of Lemma~3.2 in \cite{BC} to $S$
and $\Delta_R$ we find a compactly supported 
Hamiltonian isotopy of $\left( V, \go_{h \circ f} \right)$ 
disjoining $S$ from itself.
Theorem~2 now shows that $\gf_H$
has infinitely many nontrivial geometrically distinct
periodic points corresponding to contractible periodic orbits.  
\proofend

\section{The Weinstein conjecture}  \label{weinstein}

\ni 
Consider a weakly exact split-convex symplectic manifold $(M, \go)$.
A {\it hypersurface} $S$ in $M$
is by definition a $C^2$-smooth compact connected orientable codimension
$1$ submanifold of $M$ without boundary.
We recall that a characteristic on $S$ is an embedded circle in $S$ all
of whose tangent lines belong to the
distinguished line bundle 
\[
\cl_S \,=\, \left\{ (x, \xi) \in TS \mid \go(\xi, \eta) =0 
\text{ for all } \eta \in T_x S \right\} .
\]
We denote by $\cp^\circ (S)$ the set of closed characteristics on $S$
which are contractible in $M$.
Given $x \in \cp^\circ (S)$ we define the {\it reduced action}\, of $x$
by
\[
\ca (x) \,=\ \left| \int_{D^2} \overline{x}^* \go \right|
\]
where $\overline{x} \colon D^2 \ra M$ is a smooth disc in $M$ bounding $x$.
The {\it action spectrum}\, of $S$ is the subset 
$\gs (S) = \left\{ \ca (x) \mid x \in \cp^\circ (S) \right\}$
of $\RR$.
If $\gs (S)$ is non-empty, we define $\gl_1 (S) \in [0, \infty[$ as
\[
\gl_1 (S) \,=\, \inf \left\{ \gl \in \gs (S) \right\} .
\]
Examples show that $\gs (S)$ can be empty, see \cite{Gi, GG}.
We therefore follow \cite{HZ1} and consider parametrized neighbourhoods of $S$.
Since $S$ is orientable, there exists 
(after adding a collar $\pp M_j \times ]0,\eps]$ to each $M_j$, $j=1,
\dots, k$, in case $S$ touches $\pp M$)
an open neighbourhood $I$ of $0$ and a $C^2$-smooth diffeomorphism 
\[
\psi \colon S \times I \,\ra\, U \subset M
\]
such that $\psi (x,0) =x$ for $x \in S$.
We call $\psi$ a {\it thickening of $S$}, and we
abbreviate $S_\eps = \psi \left( S \times \left\{ \eps \right\} \right)$ and
shall often write $\left( S_\eps \right)$ instead of 
$\psi \colon S \times I \ra U$.

\begin{theorem}  \label{t:dense}
Assume that $S$ is a displaceable hypersurface of a 
weakly exact split-convex symplectic manifold $(M, \go)$,
and let $\left( S_\eps \right)$ be a thickening of $S$. 
For every $\gd > 0$ there exists $\eps \in \left[ -\gd, \gd \right]$
such that 
\[
\cp^\circ \left( S_\eps \right) \neq \emptyset 
\quad \text{and} \quad
\gl_1 \left( S_\eps \right) \le 2 e \left( S,M \right) + \gd .
\]
\end{theorem}

\proof
Fix $\gd > 0$. 
We choose $K \in \ch_c(M)$ such that $\gf_K$ displaces $S$ and $\left\|
K \right\| < e \left( S,M \right) + \gd / 2$.
Let $\rho \in \;]0,\gd]$ be so small that $\gf_K$ displaces 
the whole neighbourhood 
$\cn_\rho := \psi \left( S \times [-\rho, \rho] \right)$ of $S$.
If $\bigcup_{i \ge 1} M_i$ is an exhaustion of $M$, we choose $i$ so
large that $\supp \gf_K \subset M_i$.
We abbreviate $E = 2 e \left( S,M \right) +\gd$ and
choose a $C^\infty$-function $f \colon \RR \ra [0,E]$ such that
\[
f(t) =0 \,\text{ if } t \notin [-\rho, \rho] ,
\quad
f(0) = E,
\quad
f'(t) \neq 0 \,\text{ if } t \in \;]-\rho, \rho[ \setminus \{0\} .
\]
We define the time-independent Hamiltonian $H \in \ch_c \left( M_i
\right)$ by
\begin{equation*}  
H (x) \,=\, 
 \left\{
  \begin{array}{ll}
   f(t) & \text{if } x \in S_t , \\ [0.2em]
   0    & \text{otherwise}.
  \end{array}
 \right.
\end{equation*}
Since $\gf_H \neq \idd$ and since $\gf_H$ is supported in $\cn_\rho$,
we read off from (S1) of Theorem~\ref{t:norm} and 
from Corollary~\ref{c:ece} that
\begin{equation}  \label{e:wein:0ge}
0 \,<\, \gg_{M_i} (H) \,\le\, 2 \left\| K \right\| \,<\, E .
\end{equation}
Let $x^+ \in \cp_H$ and $x^- \in \cp_{H^-}$ be closed orbits for
which 
\[
c(H) = \ca_H \left( x^+ \right) \quad \text{ and } \quad
c \left( H^- \right) = \ca_{H^-} \left( x^- \right) .
\]
Proposition~\ref{p:yaron} applied to $H$ and $H^- = -H$ yields
\begin{eqnarray}  
\qquad c(H) \!\!&=&\!\! \ca_H \left( x^+ \right) \,=\, 
   - \int_{D^2} \left(\overline{x^+}\right)^* \go -
   \int_0^1 H \left( x^+(t) \right) dt \,\le\, 0 ,  \label{e:cH0} \\
\qquad c\left( H^- \right) \!\!&=&\!\! \ca_{H^-} \left( x^- \right) \,=\,
- \int_{D^2} 
       \left(\overline{x^-}\right)^* \go + \int_0^1 H \left( x^-(t) \right) dt
       \,\le\,  E . \label{e:cH2s}
\end{eqnarray}  
Notice that not both $x^+$ and $x^-$ are constant orbits. Indeed, if
they were, our choice of $H$ would yield 
$c(H) \in \left\{ 0,-E \right\}$ and $c \left( H^- \right) 
\in \left\{ 0,E \right\}$, and so $\gg (H) = c(H)+c \left(H^-\right)
\in \left\{ -E, 0, E \right\}$, contradicting \eqref{e:wein:0ge}.

\s
\ni
{\bf Case 1.} 
The orbit $x^+$ is not constant.
By construction of $H$ there exists $\eps \in [-\rho, \rho] \subset
[-\gd, \gd]$ such that 
$x^+ \in \cp^\circ \left( S_\eps \right)$.
The choice of $H$ and \eqref{e:cH0} yield 
$- \int_{D^2} \left( \overline{x^+} \right)^* \go \le E$. 
Assume that $- \int_{D^2} \left( \overline{x^+} \right)^* \go < -E$. 
Then \eqref{e:cH0} yields $c(H) < -E$, and so, together with
\eqref{e:cH2s}, $\gg (H) = c(H) + c \left( H^- \right) <0$,
contradicting \eqref{e:wein:0ge}.
We conclude that $\ca \left( x^+ \right) = \left| \int_{D^2} \left(
\overline{x^+} \right)^* \go \right| \le E$.

\s
\ni
{\bf Case 2.} 
The orbit $x^-$ is not constant.
Again we find $\eps \in [- \gd, \gd]$ such that $x^- \in \cp^\circ
\left( S_\eps \right)$, and arguing similarly as in Case 1 we find that
$\ca \left( x^- \right) \le E$.
The proof of Theorem~\ref{t:dense} is complete.
\proofend

A hypersurface $S$ is {\it stable}\, if there exists a thickening
$\left( S_\eps \right)$ of $S$ such
that the local flow $\psi_t$ around $S$ induced by 
$\psi \colon S \times I \ra U$ induces bundle isomorphisms
\[
T \psi_\eps \colon \cl_S \, \ra \, \cl_{S_\eps} 
\]
for every $\eps \in I$.
It then follows that $\psi_{-\eps} (x) \in \cp^\circ (S)$ for every $x
\in \cp^\circ \left( S_\eps \right)$.
Since $\psi_\eps \ra \idd$ in the $C^1$-topology as $\eps \ra 0$, and
since $\gs (S)$ is compact,
we conclude from Theorem~\ref{t:dense} the

\begin{corollary}  \label{c:stable}
Assume that $S$ is a displaceable stable hypersurface 
of a weakly exact split-convex symplectic manifold $(M, \go)$.
Then 
$\cp^\circ (S) \neq \emptyset$ and $\gl_1 (S) \le 2 e (S)$.
\end{corollary}

It is well known that every hypersurface of contact type is stable, see
\cite{HZ}, and so Theorem~3 follows from Corollary~\ref{c:stable}. 
Corollary~3 follows from Theorem~3 by using Cieliebak's result in
\cite{Ci} or by arguing as in the proof of Corollary~2 given in the
previous section.

\begin{example}
{\rm
We consider a stable hypersurface $S$ in $\left( \RR^{2n}, \go_0
\right)$.
If $S$ has diameter $\Diam (S)$, then 
$S$ is contained in a ball of radius $\Diam (S)$.
Since $e \left( B^{2n} (r) \right) = \pi r^2$, we find 
$e (S) \le \pi \Diam (S)^2$, and so
\[
\gl_1 (S) \,\le\, 2 \pi \Diam (S)^2 ,
\] 
improving the estimate in \cite{HZ1}.
}
\end{example}

\begin{remarks}\
{\rm
{\bf 1.}
Let $S$ be a stable hypersurface as in Corollary~\ref{c:stable}.
It is conceivable that the factor $2$ in the estimate $\gl_1(S) \le 2
e(S)$ can be omitted.
This is so if $S$ is a hypersurface of restricted contact type in
$\left( \RR^{2n}, \go_0 \right)$, see \cite{He}.
If $S$ bounds a convex domain $U \subset \RR^{2n}$, then 
$\gl_1(S) = c_{\HZ} (U) \le e(U) = e(S)$
where $c_{\HZ}$ is the Hofer-Zehnder capacity, \cite{HZ2}.

\s
\noindent
{\bf 2.}
Assume that $S \subset (M, \go)$ is a hypersurface of contact type and
that one of the following conditions is met.
\begin{itemize}
\item[$\bullet$]
$S$ is simply connected.
\item[$\bullet$]
$\go = d \gl$ is exact and $H^1(S;\RR) =0$.
\end{itemize}
Then $0 \notin \gs (S)$ and $\gs (S)$ is closed, cf.\ \cite{HZ1}.
Therefore, $\gl_1 (S) >0$.
\diam
}
\end{remarks}

Assume now that the hypersurface $S$ bounds, i.e., $S$ is the boundary
of a compact submanifold $B$ of $M$. 
If $M$ is simply connected, then any hypersurface $S \subset M$ bounds,
\cite{L}, and the same holds true if $H_{2n-1}(M;\ZZ) =0$;
in particular, any hypersurface of a Stein manifold of dimension at
least $4$ bounds.
In the following theorem, $\mu$ denotes the Lebesgue measure on $\RR$.

\begin{theorem}  \label{t:almost}
Assume that $(M, \go)$ is a weakly exact split-convex symplectic 
manifold  
and that $S \subset M$ is a displaceable $C^2$-hyper\-surface which bounds.
If $\left( S_\eps \right)$ with $\eps \in I$ is a displaceable 
thickening of $S$, then
\[
\mu \left\{ \eps \in I \mid \cp^\circ \left( S_\eps \right) \neq \emptyset
\right\} \,=\, \mu (I) .
\]
\end{theorem}

\proof
We can assume that $M$ is compact.
We can also assume the thickening $\left( S_\eps \right)$ to be chosen
such that for the sets $B_\eps$ bounded by $S_\eps$,
\[
B_\eps \subset B_{\eps'} \quad \, \text{if }\, \eps \le \eps '.
\]
In view of the relative monotonicity property \eqref{e:relmon} of
$c_{\HZ}^\circ$ 
the function $\eps \mapsto c_{\HZ}^\circ \left( B_\eps, M \right)$ is
then monotone increasing.
Since $S_\eps$ is displaceable, $B_\eps$ is also displaceable, and so,
according to Corollary~\ref{c:ece},
\[
c_{\HZ}^\circ \left( B_\eps, M \right) \,\le\, 
2\, e \left( B_\eps, M \right) \,<\, 
\infty \quad \text{ for all }\, \eps \in I .
\]
Theorem~\ref{t:almost} now follows from  repeating the proof of
Theorem~4 in \cite[Chapter 4]{HZ} with $C^2$-smooth instead of
$C^{\infty}$-smooth Hamiltonians and with $c_{\HZ}$ replaced by 
$c_{\HZ}^\circ$.
\proofend

\section{Closed trajectories of a charge in a magnetic field}  \label{magnetic}
 
\subsection{Proof of Theorem~4.A}

\ni
Let $(N,g)$ and $\left(T^*N, \go_{\gs} \right)$ be as in Theorem~4.A.
Since $\gs = d \ga$ is exact, 
$\go_{\gs} = -d \left( \gl+\pi^* \ga \right)$ is exact, 
and so $\left( T^*N, \go_\gs \right)$ is exact,
and since $\gs$ does not vanish, $\dim N \ge 2$, and so every energy level
$E_c = \left\{ H = c^2/2 \right\}$, $c>0$, is a $C^2$-hypersurface which
bounds.
We denote the sublevel set of $H$ by
\[
H^c \,=\, \left\{ (q,p) \in T^*N \mid H(q,p) = \tfrac 12 |p|^2 \le c
\right\} 
\]
and we define the norm of $\gs$ as
\[
\left\| \gs \right\| \,=\, \inf \left\{ \left\| \ga \right\| \mid \gs =
d \ga \right\} 
\]
where $\left\| \ga \right\| = \max_{x \in N} \left| \ga(x) \right|$.
In order to apply Theorem~\ref{t:almost} we need 

\begin{lemma}  \label{l:m:convex}
The symplectic manifold $\left(T^*N, \go_{\gs} \right)$ is convex.
Indeed, $H^c$ is of convex whenever 
$c > \frac 12 \left\| \gs \right\|^2$.
\end{lemma}

\proof
We choose a $1$-form $\ga$ on $N$ such that $d \ga = \gs$.
Under the symplectomorphism 
\[
\Phi \colon \left( T^*N, \go_\gs \right) 
                \ra \left( T^*N, \go_0 \right),
\quad
(q,p) \mapsto \left( q, p + \ga (q) \right)
\]
the Hamiltonian $H(q,p) = \frac{1}{2} |p|^2$ on 
$\left( T^*N, \go_\gs \right)$ 
corresponds to the Hamiltonian $H_\ga (q,p) = \frac{1}{2} |p-\ga|^2$ on 
$\left( T^*N, \go_0 \right)$.
If $c > \frac 12 \left\| \ga \right\|^2$, then the sublevel set 
$H^c_\ga = \left\{ (q,p) \mid H_\ga (q,p) \le c \right\}$ 
contains $N$, and so the Liouville vector field 
$\sum_i p_i \frac{\pp}{\pp p_i}$
for $\go_0$ intersects the boundary of $H_c^\ga$ transversally. 
Therefore, $H_c^\ga$ is convex.
It follows that $H_c = \Phi^{-1} \left( H_c^\ga \right)$ is
convex whenever $c > \frac 12 \left\| \ga \right\|^2$.
Since this is true for any $\ga$ with $d \ga = \gs$, the lemma follows.
\proofend

\begin{remark}  \label{r:notcontact}
{
\rm
Combining the identity \eqref{c=s} below with arguments from \cite{CMP} one
can show that if $N$ is orientable and different from the $2$-torus,
then $E_c$ is not of contact type if
$c \le \frac 12 \left\| \gs \right\|^2$, and so $H^c$ is not convex if
$c \le \frac 12 \left\| \gs \right\|^2$.
}
\end{remark}

Let $\chi (N)$ be the Euler characteristic of $N$.

\s
\ni
{\bf Case 1.} $\chi (N) =0$.
We set 
\begin{equation}  \label{def:dgs}
d \,=\, 
d(g, \gs) \,=\, \sup \left\{ c \ge 0 \mid H^c \text{ is displaceable in }
\left( T^*N, \go_\gs \right) \right\} .
\end{equation}
Notice that since $\dim N \ge 2$, 
\[
d \,=\, \sup \left\{ c \ge 0 \mid E_c \text{ is displaceable in }
\left( T^*N, \go_\gs \right) \right\} .
\]
Since $\gs \neq 0$, the zero section $N$ of $T^*N$ is not Lagrangian,
and so a remarkable theorem of Polterovich \cite{P2,LS} implies that $d>0$.
We shall see below that $d<\infty$.
Theorem~4.A follows from applying Theorem~\ref{t:almost} to $S = E_{d/2}$
and a thickening
\[
\psi \colon S \times \left] -d/2,d/2 \right[  \,\ra\, \bigcup_{0<c<d} E_c
\]
such that $\psi \left( S \times \{ \eps \} \right) = E_{\eps + d/2}$.

\s
\ni
{\bf Case 2.} $\chi (N) \neq 0$.
In this case the zero section $N$ is not displaceable for topological
reasons. 
We use a stabilization trick used before by Macarini \cite{Mac}.
Let $S^1$ be the unit circle, and  denote canonical coordinates on
$T^*S^1$ by $(x,y)$.  
We consider the manifold 
$T^* \left( N \times S^1 \right) = T^* N \times T^* S^1$ 
endowed with the split symplectic form 
$\go = \go_\sigma \oplus \go_{S^1}$, 
where $\go_{S^1} = dx \wedge dy$.
In view of Lemma~\ref{l:m:convex}, 
$\left( T^*N \times T^* S^1, \go \right)$
is a weakly exact convex symplectic manifold.
Moreover, $N \times S^1$ is not Lagrangian, and $\chi \left( N \times
S^1 \right) =0$.
Let
\[
H_1(q,p) = \tfrac{1}{2} |p|^2,
\quad
H_2(x,y) = \tfrac{1}{2} |y|^2,
\quad
H(q,p,x,y) =  \tfrac 12 |p|^2 + \tfrac 12 |y|^2 
\]
be the metric Hamiltonians on $T^*N$, $T^*S^1$ and $T^*N \times T^*S^1$.
In order to avoid confusion, we denote their energy levels by
$E_c(H_1)$, $E_c(H_2)$ and $E_c(H)$.
Repeating the argument given in Case~1 for the Hamiltonian system
\begin{equation}  \label{e:m:H}
H \colon \left( T^*N \times T^* S^1, \go \right) \ra \RR 
\end{equation} 
and 
\[
d \,=\, d(g, \gs) \,=\, \sup \left\{ c \ge 0 \mid H^c \text{ is
displaceable in } \left( T^*N \times T^* S^1, \go \right) \right\} 
\]
we find that 
\[
\mu \left\{ \eps \in \;]0, d[ \mid \cp^\circ \left( E_\eps(H) \right)
\neq \emptyset \right\} \,=\, d .
\]
Fix $\eps \in \;]0,d[$ such that $\cp^\circ \left(E_\eps (H) \right)
\neq \emptyset$.
Since the Hamiltonian system \eqref{e:m:H} splits,
a contractible closed orbit $x(t)$ on $E_\eps (H)$ is of the form
$\left( x_1(t), x_2(t) \right)$, where $x_1$ is a contractible
closed orbit on $E_{\eps_1} (H_1)$ and $x_2$ is a contractible closed orbit
on $E_{\eps_2} (H_2)$ and $\eps_1 + \eps_2 = \eps$. 
Since the only contractible orbits of $H_2 \colon T^*S^1 \ra \RR$
are the constant orbits on $E_0 (H_2)$, we conclude that $\eps_2 = 0$ 
and $\eps_1 = \eps$, and so $x_1 \in \cp^\circ \left( E_\eps (H_1)\right)$.
It follows that
\[
\mu \left\{ \eps \in \;]0, d[ \mid \cp^\circ \left( E_\eps (H_1)\right)
\neq \emptyset \right\} \,=\, d .
\]
The proof of Theorem~4.A is complete.
\proofend

\subsection{Comparison of $d(g,\gs)$ and $\frac 12 \left\| \gs \right\|^2$}

It would be important to know a computable lower bound of $d(g, \gs)$.
An upper bound can be described in a variety of ways.

\begin{proposition}  \label{p:mane}
We have
$d (g,\gs) \le \frac 12 \left\| \gs \right\|^2$.
\end{proposition}

\proof
We assume first that $\chi (N) =0$. 
Arguing by contradiction, we assume that 
$d = d(g,\gs) > \frac 12 \left\| \gs \right\|^2$.
We then find a $1$-form $\ga$ on $N$ such that $d \ga = \gs$ and $d >
\frac 12 \left\| \ga \right\|^2$.
By definition of $d$, the graph $\Gamma_{-\ga}$ of $-\ga$, which is
contained in $H^{\frac 12 \| \ga \|^2}$, is then a displaceable subset of 
$\left( T^*N, \go_{\gs} \right)$, 
and so the zero section $\Phi \left( \Gamma_{-\ga} \right)$
of $T^*N$ is a displaceable subset of $\left( T^*N, d\gl \right)$.
This contradicts a Lagrangian intersection result of Gromov \cite{Gr}.

Assume now that $\chi (N) \neq 0$. We denote by $g_{S^1}$ the Riemannian
metric of the unit circle. By definition of $d (g, \gs)$ and by the
already proved case,
\[
d(g,\gs) \,=\, d \left( g \oplus g_{S^1} , \gs \oplus 0 \right)
         \,\le\, \tfrac 12 \left\| \gs \oplus 0 \right\|^2 
         \,\le\, \tfrac 12 \| \gs \|^2 .
\]
The proof of Proposition~\ref{p:mane} is complete.
\proofend

An important number associated with the Hamiltonian system~\eqref{e:Hmag}
is 
{\it Ma\~n\'e's strict critical value} $c_0(g,\gs)$ for whose definition
and relevance we refer to \cite{PP,CIPP1,PPS}.
Let $\ga$ be such that $d \ga = \gs$. 
According to Corollary~1 in \cite{CIPP1}, $c_0(g, \gs)$ is given by
\begin{equation}  \label{e:mane}
c_0(g, \gs) \,=\, \inf \max_{x \in N} \tfrac 12 \left| \gb - \ga \right|^2
\end{equation} 
where the infimum is taken over all closed $1$-forms $\gb$ on $N$.
It follows that
\begin{equation}  \label{c=s}
c_0 (g,\gs) \,=\, \tfrac 12 \left\| \gs \right\|^2 .
\end{equation}
We denote by $\Lambda_{-\ga}$ the set of Lagrangian submanifolds in
$\left( T^*N, \go_\gs \right)$ which are Lagrangian isotopic to the graph
$\Gamma_{-\ga}$ of $-\ga$. Combining \eqref{e:mane} with a result in
\cite{PPS}, we find
\[
c_0(g,\gs) \,=\, \inf \left\{ c \in \RR \mid H^c \text{ contains a
Lagrangian submanifold in } \Lambda_{-\ga} \right\} .
\]
This is a purely symplectic characterization of 
$c_0 (g, \gs) = \frac 12 \left\| \gs \right\|^2$.

\s
We recall from Theorem~4.A that $\cp^\circ \left(E_c \right) \neq
\emptyset$ for almost all $c \in \; \left] 0, d(g, \gs) \right]$.
It follows from Lemma~\ref{l:m:convex} and a theorem of Hofer and
Viterbo \cite{HV} that $E_c$ carries a closed orbit whenever 
$c > \frac 12 \| \gs \|^2$.
More precisely, for every non-trivial homotopy class $h \in \pi_1(N)$ and
every $c > c_0(g, \gs) = \frac 12 \left\| \gs \right\|^2$ there exists a
closed orbit on $E_c$ whose projection to $N$ lies in $h$, see
\cite[Theorem 27]{CIPP2}.
The following example shows that $\cp^\circ (E_c)$ can be empty for all
$c \ge \frac 12 \left\| \gs \right\|^2$.
It also shows that there can be a gap between $d(g,\gs)$ and $\frac 12
\left\| \gs \right\|^2$.

\begin{example}  \label{ex:12}
{\rm
Let $N$ be a closed orientable surface of genus $2$.
It has been shown in \cite{PP} that there exists a Riemannian metric $g$
and an exact $2$-form $\gs$ on $N$ such that
\begin{itemize}
\item[(i)]
$c_0(g,\gs) > \frac 12$;
\item[(ii)]
the restriction of the flow of \eqref{e:Hmag} to $E_c$ is Anosov for all
$c \ge \frac 12$.
\end{itemize}
Property (ii) implies that $\cp^\circ (E_c) = \emptyset$ for all $c \ge
\frac 12$, and so, by Theorem~4.A, Property~(i) and \eqref{c=s},
\[
d (g, \gs) \,\le\, \tfrac 12 \,<\, 
c_0(g,\gs)  \,=\, \tfrac 12 \left\| \gs \right\|^2 .
\]
}
\end{example}

\subsection{Proof of Theorem~4.B}  \label{ss:magnetic:4B}

\ni
We say that a closed $2$-form $\gs$ on a manifold $N$ is {\it
rational}\,
of
\[
\hbar \,:=\, \inf_{\left[ S \right] \in \pi_2 (N)} 
\left\{ \int_S \gs \;\Bigg|\; \int_S \gs >0 \right\} \,>\, 0 .
\]
Our most general result about the existence of closed orbits of magnetic
flows is
\begin{theorem}  \label{t:maggen}
Assume that $N = N_1 \times N_2 \times N_3$ is a closed manifold, where
$N_1$ is any closed manifold, $N_2 = \times_i S^2$ is a product of
$2$-spheres, and $N_3 = \times_j \Sigma_j$ is a product of closed
orientable surfaces of genus at least $2$,
and assume that $N$ is endowed with a $C^2$-smooth 
Riemannian metric $g$ and a non-vanishing closed $2$-form $\gs$ such
that 
\[
[\gs] \,=\, 0 \oplus \left[ \gs_2 \right] \oplus \left[ \gs_3 \right]
\,\in\, H^2 \left( N_1 \times N_2 \times N_3 \right) ,
\] 
such that $\left[ \gs_2 \right]$ is rational, 
and such that
$\left[ \gs_3 \right] \in H^2 \left( N_3 \right)$ is cohomologically
split in the sense that
\[
\left[ \gs_3 \right] \,\in\, \oplus_i \RR \left[ \Sigma_i \right] \,=\, 
\oplus_i H^2 \left( \Sigma_i \right) 
\,\subset\, H^2 \left( \displaystyle\times_i \Sigma_i \right) .
\]
\begin{itemize}
\item[(i)]
If $\left[ \gs_2 \right] \neq 0$, there exists $d>0$ such that $\cp^\circ
\left(E_c\right) \neq \emptyset$ for a dense set of values $c \in
\;]0,d]$.
\item[(ii)]
If $\left[ \gs_2 \right] = 0$, there exists $d>0$ such that 
$\cp^\circ \left(E_c\right) \neq \emptyset$ for almost all $c \in \;]0,d]$.
\end{itemize}

\end{theorem}
\ni
For $N_2$ and $N_3$ a point, Theorem~\ref{t:maggen} is Theorem~4.A, and
for $N_1$ a point and $N_3$ or $N_2$ a point, Theorem~\ref{t:maggen} is
a generalization of Theorem~4.B\:(i) or (ii).

\m
\ni
{\it Proof of Theorem~\ref{t:maggen}:}
We first consider a closed orientable surface $\Sigma$ different from
the torus, and we endow $\Sigma$ with a Riemannian metric $g$ of
constant curvature $k$.
We fix an orientation of $\Sigma$, denote the area form on $\Sigma$ by
$\tau$, and consider the $2$-form $\gs = s \tau$ for some $s \in \RR$.
Recall that $\go_\gs = \go_0 - \pi^* \gs$.
The following lemma was explained to us by Viktor Ginzburg.

\begin{lemma}  \label{l:sigma:convex}
The symplectic manifold $\left( T^*\Sigma, \go_\gs \right)$ is convex.
Indeed, if $\Sigma = S^2$, then $H^c$ is convex for all $c>0$, 
and if $\genus ( \Sigma ) \ge 2$, then $H^c$ is convex for all 
$c > -\frac{s^2}{2k}$.
\end{lemma}

\proof
We fix $c>0$ and consider $E_c$ as an oriented $S^1$-bundle 
\begin{equation}  \label{e:SEN}
S^1 \,\longrightarrow\, E_c \,\stackrel{\pi_c}{\longrightarrow}\, N .
\end{equation}
Let $X_c$ be the geodesic spray on $E_c$,
let $Y_c$ be the vector field on $E_c$ generating the $S^1$-action,
and let $\ga_c$ be the connection $1$-form of the bundle \eqref{e:SEN}.
Then
\begin{equation}  \label{e:connection}
\ga_c \left( X_c \right) = 0, \quad
\ga_c \left( Y_c \right) = 1, \quad
d \ga_c = - \pi_c^* (k \tau) .
\end{equation}
Varying over $c>0$ we obtain vector fields $X,Y$ 
and a $1$-form $\ga$ on $T^*N \setminus N$ such that $\ga |_{E_c} =
\ga_c$ and $d \ga = - \pi^* ( k \tau )$.
Since $N$ is not the torus, $k \neq 0$, and so we can set $\gb =
-\frac{s}{k} \ga$.
Then 
\[
d \gb \,=\, - \tfrac{s}{k} d \ga \,=\, \pi^* (s \tau) \,=\,
\pi^* \gs \quad\, \text{on }\, T^*N \setminus N .
\]
Therefore,
\begin{equation}  \label{e:dlb}
d \left( - \gl - \gb \right) \,=\, \go_\gs .
\end{equation}
The vector field $X_H = X-sY$ on $T^*N \setminus N$ is the Hamiltonian
vector field of $H(q,p) = \frac 12 |p|^2$ with respect to $\go_\gs$.
In particular, $X_H |_{E_c}$ is a section of the distinguished line
bundle $\cl_{E_c}$ for every $c>0$.
Notice that
$\gl (X) |_{E_c} = 2c$ and $\gl (Y) =0$.
Moreover, $\gb = - \frac sk \ga$ and \eqref{e:connection} yield 
$\gb (X) = 0$ and $\gb (Y) = - \frac sk$.
Therefore,
\begin{equation}  \label{lbX}
\left( -\gl - \gb \right) \left( X_H \right) \,=\, -2c - \tfrac{s^2}{k} .
\end{equation}
Equation \eqref{e:dlb} and \eqref{lbX} show that if $N=S^2$, then $E_c$
if of contact type for every $c>0$, 
and if $\genus (N) \ge 2$, then $E_c$ is of contact type if 
$c \neq \frac{s^2}{2k}$.
If $s=0$, all these hypersurfaces are convex, and so the claim follows.
\proofend

Let now $N$, $g$ and $\gs$ be as in Theorem~\ref{t:maggen}.
We denote the area form $\tau$ considered in Lemma~\ref{l:sigma:convex}
by $\tau_{S^2}$ or $\tau_{\Sigma}$.
By assumption on the form $\gs_3$ there are real numbers $s_i$ and $s_j$
such that
\[
\left[ \gs_2 \right] \,=\, \oplus_i s_i \left[ \tau_{S^2} \right]
\,\in\,
H^2 \left( \times_i S^2 \right),
\quad\,
\left[ \gs_3 \right] \,=\, \oplus_j s_j \left[ \tau_{\Sigma_j} \right]
\,\in\,
H^2 \left( \times_j \Sigma_j \right) .
\]
Define the closed $2$-form $\gs_0$ on $N = N_1 \times N_2 \times N_3$ as
\[
\gs_0 \,=\, 0 \oplus_i s_i \tau_{S^2} \oplus_j s_j \tau_{\Sigma_j} . 
\]
According to Lemma~\ref{l:sigma:convex} the symplectic manifold
\[
\left( T^*N, \go_{\gs_0} \right) \,=\,
\left( T^* N_1, \go_0 \right) \times_i \left( T^*S^2, s_i \tau_{S^2}
\right) \times_j \left( T^* \Sigma_j, s_j \tau_{\Sigma_j} \right)
\]
is a product of convex symplectic manifolds.
By assumption on $\gs$ there exists a $1$-form $\ga$ on $N$ such that
$\gs \,=\, \gs_0 + d\ga$.
The next lemma will allow us to interpolate between the forms $\go_\gs$
and $\go_{\gs_0}$.

\begin{lemma}  \label{l:inter}
For every $r>0$ there exists $R>0$ and a smooth function $f \colon \RR
\ra [0,1]$ such that
\begin{equation}  \label{es:f}
f(t) = 1, \,\, s \le r; 
\quad\,
f(t) = 0, \,\, s \ge R ,
\end{equation}
and such that the closed $2$-form $\go_f$ defined as
\[
\go_f (q,p) \,:=\, \go_{\gs_0}(q,p) - 
d \,\big( f \left( |p| \right) \pi^* \ga (q) \big)
\]
is nondegenerate and hence symplectic on $T^*N$.
\end{lemma}

\proof 
Fix $(q,p) \in T^*N$.
For convenience we choose local coordinates $q_i$ around $q$ on $N$ such
that for the coefficients $g_{ij}$ of $g$ we have $g_{ij}(q) = \gd_{ij}$
and $g_{ij,k}(q) = 0$ for all $i,j,k$.
Let $\gs_0$ and $\ga$ be given by $\gs_0 (q) = \sum_{i,j} S_{ij}(q) \,dq_i
\wedge dq_j$ and $\ga (q) = \sum_i A_i(q) \,dq_i$. Then
\begin{eqnarray*}
\go_f (q,p) &=&
\sum_{i,j} \left( \gd_{ij} + A_i(q) f' \left( |p| \right)
\tfrac{p_j}{|p|} \right) dq_i \wedge dp_j  \\
 & & + \sum_{i,j} \left( A_{i,j} (q) f \left( |p| \right)  - S_{ij}(q)
 \right) dq_i \wedge dq_j.  
\end{eqnarray*}
The square root of the determinant of the matrix of $\go_f (q,p)$, which
we want to be non-zero, is therefore
\begin{equation}  \label{e:detda}
\det \left( \gd_{ij} + A_i(q) f'\left( |p| \right) \tfrac{p_j}{|p|} \right).
\end{equation}
Choose $\eps >0$ so small that $\det \left( \gd_{ij} + c_{ij} \right)
>0$ whenever $\left| c_{ij} \right| \le \eps$ for all $i,j$.
Since $N$ is compact, we find $a < \infty$ such that for every $q \in N$
there exists a Riemannian metric as above such that $\left| A_i(q)
\right| \le a$ for all $i$.
Choose now $R> r + a/\eps$ and $f \colon \RR \ra [0,1]$ satisfying
\eqref{es:f} and $\left| f'(r) \right| \le \eps / a$.
Then
$\left| A_i(q) f' \left( |p| \right) \frac{p_j}{|p|} \right| \le \eps$,
and so the determinant~\eqref{e:detda} does not vanish.     
\proofend

\s
\ni
{\bf Case 1.} $\chi (N) = 0$.
Since $\gs \neq 0$, the full result of \cite{P2,LS} implies that
the displacement energy of $N$ in $\left( T^*N, \go_\gs \right)$
vanishes.
We therefore find $d >0$ such that 
$e \left( H^d, T^*N \right) \le \hbar /2$. 
Fix $d' \in \;]0,d[$ and choose 
$\gf \in \Ham_c \left( T^*N, \go_\gs \right)$ displacing $H^{d'}$.
Choose $r>0$ so large that $\supp \gf \subset T_r^*N$,
and then choose $R$ and $f$ as in Lemma~\ref{l:inter}.
With these choices,
$\gf \in \Ham_c \left( T^*N, \go_f \right)$, and $\gf$
displaces $H^{d'}$ in $\left( T^*N, \go_f \right)$.
Moreover, $\go_f = \go_{\gs_0}$ on $T^*N \setminus T_R^* N$, and so
$\left( T^*N, \go_f \right)$ is split-convex.

\m
(i) If $\left[ \gs_2 \right] \neq 0$, then $\left( T^*N, \go_f \right)$
is not weakly exact.
However, we have

\begin{lemma}
The first Chern class $c_1 \left( T^*N, \go_f \right)$ vanishes on
$\pi_2 \left( T^*N \right)$.
\end{lemma}
\proof
We abbreviate $M = T^* N$.
The tangent bundle of $M$ at a point $(q,0) \in N$
naturally splits as $T_{(q,0)}M \cong T^*_q N \oplus T_q N$. 
Notice that the summand $T^*N$ is a Lagrangian subbundle 
of the restriction of $TM$ to $N$ for the symplectic structure 
$\go_f = \go_\gs = -d \gl - \pi^* \gs$.
Therefore, $c_1 \left( M, \go_f \right)$ vanishes on 
$\pi_2 (N) = \pi_2 \left( M \right)$. 
\proofend

Notice that $\hbar \left( \go_f \right) = \hbar \left( \go_\gs \right) =
\hbar (\gs_2)$.
According to Theorem~\ref{t:app:dense},
$\cp^\circ \left( E_c \right) \neq \emptyset$ for a dense set of $c \in
\;]0,d']$. Since $H^{d'} \subset \supp \gf \subset T_r^*N$, we have
$\go_f = \go_\gs$ on $H^{d'}$, and so these closed characteristics are
characteristics with respect to the original symplectic structure
$\go_\gs$.
Since $d' \in \;]0,d[$ was arbitrary, Theorem~\ref{t:maggen}\:(i) 
for $\chi (N)=0$ follows.

\s
(ii) If $\left[ \gs_2 \right] =0$, then $\left( T^*N, \go_f \right)$ is
a weakly exact split-convex symplectic manifold.
Applying Theorem~\ref{t:almost} and using that 
$d' \in \;]0,d[$ was arbitrary, Theorem~\ref{t:maggen}\:(ii) for
$\chi (N)=0$ follows.

\s
\ni
{\bf Case 2.} $\chi (N) \neq 0$. We can now take 
$d = d (g,\gs)$ as in \eqref{def:dgs}.
We stabilize $\left( T^*N, \go_\gs \right)$ by $\left( T^* S^1,
\go_{S^1} \right)$ as in Case~2 of the proof of Theorem~4.A
and combine the arguments there with the arguments in Case~1 above.
The proof of Theorem~\ref{t:maggen} is complete. 
\proofend

\begin{remarks}  \label{r:dfinite}
{\rm
In view of Example~\ref{ex:12}, 
the number $d>0$ in Theorem~4.B (ii) cannot be chosen arbitrarily large
in general. 
Here is a simpler example illustrating this fact:
Let $N$ be a closed oriented surface equipped with a metric of constant
curvature $-1$, and let $\gs$ be the area form on $N$.
If $c \ge \frac 12$, then $\cp^\circ \left( E_c \right) = \emptyset$,
see \cite[Example 3.7]{Gi0}.
}
\end{remarks}

\subsection{The state of the art}  \label{ss:state}

\s
\ni
We shall only consider the existence problem of closed orbits on small
energy levels and refer to the review \cite{Gi0} for results concerned
with closed orbits on intermediate and large energy levels.

\m
\ni
{\bf 1. $\gs$ is exact.}
Theorem~4.A improves a result of Polterovich and Macarini 
\cite{P0,Mac} 
who proved $\cp^\circ (E_c) \neq \emptyset$ for a sequence $c \ra 0$.

\m
\ni
{\bf 2. $\gs$ is neither exact nor symplectic.}
The only previous results for such magnetic fields are the 
Polterovich--Macarini result stating that if $[\gs] |_{\pi_2(N)} =0$,
then $\cp^\circ (E_c) \neq \emptyset$ for a sequence $c \ra 0$,
and a result of Lu \cite{Lu} stating that for the torus $T^n$
endowed with any Riemannian metric, 
$\cp^\circ (E_c) \neq \emptyset$ for almost all $c> 0$.
Theorem~\ref{t:maggen} is thus new.

\m
\ni
{\bf 3. $\gs$ is symplectic.}
Most previous results where obtained for symplectic forms.
We refer to \cite{Gi0, GK2, GG} for the best known results
and only mention two of them.

\s
(i) If $N$ is a surface, then for {\it all}\, sufficiently small $c>0$
the energy level $E_c$ carries a closed orbit, and if $N$ is a sphere or
a torus, these orbits can be chosen in $\cp^\circ \left( E_c \right)$,
see \cite{Gi0}.

\s
(ii) If $[\gs] |_{\pi_2 (N)} =0$, then $\cp^\circ (E_c) \neq \emptyset$
for almost all sufficiently small $c>0$, see \cite{GG}.

\m
\ni
Theorem~\ref{t:maggen} is new if $\left[ \gs_2 \right] \neq 0$ and $\dim
N \ge 4$.

\s
\b
The $2$-sphere and the $2$-torus are of particular interest, see \cite{No,Ko}.

\begin{corollary}
Assume that $g$ is a $C^2$-smooth Riemannian metric on $S^2$ and that
$\gs \neq 0$ is a closed $2$-form on $S^2$.
Then $\cp^\circ (E_c) \neq \emptyset$ for a dense subset of small $c>0$
and for all sufficiently small $c>0$ if $\gs$ is symplectic.
\end{corollary}

\begin{corollary}
Assume that $g$ is a $C^2$-smooth Riemannian metric on $T^2$ and that
$\gs \neq 0$ is a closed $2$-form on $T^2$.
Then $\cp^\circ (E_c) \neq \emptyset$ for almost all sufficiently small
$c>0$ and for all sufficiently small $c>0$ if $\gs$ is symplectic.
\end{corollary}

\begin{example}
{\rm
We consider a non-vanishing magnetic potential $A$ on Euclidean space 
$\RR^3$ which is parallel to the $z$-axis 
and is $2\pi$-periodic in $x$ and $y$.
The induced closed $2$-form $\gs = A(x,y) \, dx \wedge dy$ on the flat
torus $T^2 = \left\{ x,y \!\mod 2\pi \right\}$ is exact if and only if
\[
\int_{T^2} A(x,y) \,dx dy \,=\, 0 
\]
and symplectic if and only if $A(x,y) \neq 0$ for all $(x,y) \in T^2$.
If $\gs$ is exact, $\cp^\circ (E_c) \neq \emptyset$ for
almost all $c \in \;]0,d (T^2, \gs)]$ by Theorem~4.A,
if $\gs$ is neither exact nor symplectic,  
$\cp^\circ (E_c) \neq \emptyset$ for almost all $c>0$ 
by a result of Lu \cite{Lu}, 
and if $\gs$ is symplectic,
$\cp^\circ (E_c) \neq \emptyset$ for all $c>0$
by a result of Arnold applying the Conley--Zehnder theorem,
see \cite[Theorem 3.1\:(i)]{Gi0}.
The projections of all these closed trajectories lift to closed
trajectories of speed $c$ in the $x$-$y$-plane of a charge subject
to the magnetic potential $A$. 
}
\end{example}

\section{Lagrangian intersections}  \label{lag}

\ni
Theorem~5 is a special case of 
\begin{theorem}  \label{t:lag:prod}
Assume that $(M,\go)$ is a product of weakly exact convex symplectic
manifolds, 
and let $L \subset M \setminus \pp M$ be a closed Lagrangian submanifold
such that
\begin{itemize}
\item[(i)]
the injection $L \subset M$ induces an injection $\pi_1(L) \subset
\pi_1(M)$;  
\s
\item[(ii)]
$L$ admits a Riemannian metric none of whose closed geodesics is contractible.
\end{itemize}
Then $L$ is not displaceable.
\end{theorem}

\proof
Arguing by contradiction we assume that $\psi \in \Ham_c(M,\go)$ 
displaces $L$.
We can assume that $M$ is compact. 
By Weinstein's Theorem we find $\eps >0$ such that a neighbourhood
$U_{\eps}$ of
$L$ in $M$ can be symplectically identified with $T_{3\eps}^*L$.
Choose a smooth function $f \colon [0,3\eps] \ra [0,1]$ such that
\[
f(r) = -1 \,\text{ if } r \le \eps, 
\quad
f(r) = 0 \,\text{ if } r \ge 2 \eps, 
\quad
f'(r) > 0 \,\text{ if } r \in \;]\eps, 2 \eps[ . 
\]
We choose canonical coordinates $(q,p)$ on 
$T^*_{3\eps} L \equiv U_{\eps}$ and define the
autonomous Hamiltonian $H \colon M \ra \RR$ by
\[
H(x) = H(p) = f \left( |p| \right) \,\text{ if } x = (q,p) \in U_{\eps},
\quad
H(x) = 0 \,\text{ otherwise}.
\]
Set again $H^{(n)} (t,x) = n H (nt,x)$ so that $\gf_{H^{(n)}} = \gf_H^n$.
By assumptions (i) and (ii) and by our choice of $H$, 
the only contractible periodic orbits of
$\gf_H^t$ are fixed points, and so
$\Sigma_{H^{(n)}} = \left\{ 0,n \right\}$.
Since $\gf_H^n \neq \id$, $\gg \left( \gf^n_H \right) >0$, and so we
conclude that
\begin{equation}  \label{e:infty}
\gg \left( \gf_H^n \right) =n \ra \infty \quad \text{as } n \ra \infty .
\end{equation}
We now choose $\eps >0$ above so small that 
$\psi \left( U_{\eps} \right) \cap U_{\eps} = \emptyset$.
Since $\gf_H^n$ is supported in $U_{\eps}$ for all $n$, we conclude from
Proposition~\ref{p:phi2psi} that $\gg \left( \gf_H^n \right) \le 2 \gg (\psi)$, which by
\eqref{e:infty} is a contradiction.
\proofend

\begin{remarks}
{\rm
{\bf 1.}
The conclusion of Theorem~\ref{t:lag:prod} does not hold for a small circle
$L$ in a disc $D^2$, showing that condition~(i) cannot be omitted.

\s
{\bf 2.}
According to a theorem of Gromov, \cite[$2.3.\text{B}_3'$]{Gr},
the conclusion of Theorem~\ref{t:lag:prod} holds for any closed Lagrangian
submanifold $L \subset M \setminus \pp M$ for which 
$[\go] |_{\pi_2(M,L)}=0$.

\s
{\bf 3.}
In \cite{LP}, Lalonde and Polterovich used the general energy-capacity
inequality to prove the conclusion of Theorem~\ref{t:lag:prod} for
{\it any}\, symplectic manifold $(M, \go)$ and any closed Lagrangian
submanifold $L \subset M \setminus \pp M$ satisfying (i) and
\begin{itemize}
\item[(ii')]
$L$ admits a Riemannian metric of non-positive curvature.
\end{itemize}
Of course, (ii') implies (ii).
We show by an example that (ii) is a weaker condition than (ii').
Let $H$ be the $(2k+1)$-dimensional Heisenberg group endowed with any left
invariant Riemannian metric, and 
choose a discrete cocompact subgroup $\Gamma \subset H$.
The Riemannian exponential map from the Lie algebra of $H$ to $H$ is not
injective, but there are no closed geodesics, see e.g.\ \cite{E}.
Therefore, $\Gamma \setminus H$ satisfies condition (ii).
On the other hand, $\pi_1(\Gamma \setminus H) = \Gamma$ is nilpotent, and so 
$\Gamma \setminus H$ cannot satisfy (ii'), see \cite{GW,Y}.
}
\end{remarks}

\appendix \ 

\section{An extension to semi-positive convex symplectic manifolds}

\ni
In this appendix we extend parts of the main body of this paper to more
general convex symplectic manifolds than weakly exact ones.
We shall not aim at outermost generality but shall focus on those
additional results needed for the proof of
Theorem~\ref{t:maggen}.

\s
For any symplectic manifold $(M, \go)$ with first Chern class $c_1 = c_1
(\go)$ we denote the homomorphisms $\pi_2(M) \ra \RR$ defined by
integration of $\go$ and a representative of $c_1$ over a sphere also by
$\go$ and $c_1$.
Following \cite{M1,HS} we say that a $2n$-dimensional symplectic manifold 
$(M,\go)$ is {\it semi-positive}\, if one of the following conditions 
is satisfied.
{\rm
\begin{itemize}
\item[(SP1)] 
  $\go (A) = \gl\, c_1 (A)$ for every $A \in \pi_2(M)$ where $\gl \ge 0$;
\s
\item[(SP2)] 
  $c_1(A)=0$ for every $A \in \pi_2(M)$;
\s
\item[(SP3)] 
  The minimal Chern number $N \ge 0$ defined by 
  $c_1 \left( \pi_2(M) \right) = N \ZZ$ is at least $n-2$.
\end{itemize}
The semi-positivity condition will exclude bubbling off of pseudo
holomorphic spheres in the compactifications
of the moduli spaces relevant for defining Floer homology.

Consider a semi-positive compact split-convex symplectic manifold $(M,
\go)$.
We denote by $\cl$ the space of smooth contractible loops $x \colon S^1
\ra M$.
If $(M, \go)$ is not weakly exact, then the action functional
\eqref{def:af} is not well-defined on $\cl$.
The action functional is, however, well-defined on a suitable cover
$\widetilde{\cl}$ of $\cl$.
The elements of $\widetilde{\cl}$ are 
equivalence classes $[x,\bar{x}]$ of pairs $(x,\bar{x})$ 
where $x \in \cl$ and 
$\bar{x} \colon D^2 = \{z \in \CC \mid |z| \le 1\} \ra M$ 
satisfies $\bar{x}(e^{it})=x(t)$,
and where $(x_1,\bar{x}_1)$ and $(x_2,\bar{x}_2)$ are equivalent if 
\[
x_1 = x_2, \quad 
\go (\bar{x}_1 \# \bar{x}_2) =0, 
\quad
c_1 (\bar{x}_1 \# \bar{x}_2) =0.
\]
The group
\[
\Gamma \,=\, \frac{\pi_2(M)}{\ker(c_1)\cap \ker(\omega)}
\]
acts on pairs $[x,\bar{x}]$ by
\[
[x,\bar{x}] \,\mapsto\, [x,\bar{x}\#A], \quad\, A \in \Gamma ,
\]
and $\cl = \widetilde{\cl} / \Gamma$. 
The action functional
\[
\widetilde{\ca}_H([x,\bar{x}]) \,:=\, -\int \bar{x}^*\omega-\int_0^1
H_t(x(t)) \,dt
\]
is well-defined on $\widetilde{\cl}$. 
For an admissible Hamiltonian function $H \in \ch$ the set of its 
critical points $\widetilde{\cp}_H$ consists of those 
$[x,\bar{x}] \in \widetilde{\cl}$
for which $x \in \cp_H$ is a $1$-periodic orbit of the flow $\gf_H^t$. 
Even for regular admissible Hamiltonians $H \in \ch_{\reg}$
the action functional $\widetilde{\ca}_H$ will have infinitely many
critical points, and so we need to define the Floer chain complex
over a Novikov ring. 
The Novikov ring $\Lambda_\Gamma$ consists of finite formal sums
\[
\sum_{\gg \in \Gamma}r_\gg \,\gg, \quad\, r_\gg \in \ZZ_2 ,
\]
which satisfy the finiteness condition
\[
\# \left\{ \gg \in \Gamma \mid r_\gg \neq 0, \,\,
\omega(\gamma)\geq \kappa \right\} 
\,<\, \infty \quad\, \text{for all }\, \kappa \in \RR.
\]
The Novikov ring $\Lambda_\Gamma$ is naturally graded by $-2c_1$. 
Since its coefficients are taken in the field $\ZZ_2$, 
the Novikov ring $\Lambda_\Gamma$ is actually a field. 
For $H \in \ch_{\reg}$ we define
$CF (M;H)$ to be the $\Lambda_\Gamma$-vector
space consisting of formal sums
\[
\sum_{[x,\bar{x}] \in \widetilde{\cp}_H}r_{[x,\bar{x}]} \,[x,\bar{x}],
\quad r_{[x,\bar{x}]} \in \ZZ_2 ,
\]
which meet 
\[ 
\# \left\{ [x,\bar{x}] \in \widetilde{\cp}_H \;\big|\;  r_{[x,\bar{x}]} \neq 0,
\,\
\widetilde{\ca}_H \left( [x,\bar{x}] \right) \le \kappa \right\} \,<\, \infty 
\quad\,
\text{for all }\, \kappa \in \RR .
\]
For $[x,\bar{x}] \in \widetilde{\cp}_H$
there is a well defined Conley--Zehnder index $\mu$, which satisfies 
\[
\mu([x,\bar{x}]\#A) \,=\, \mu([x,\bar{x}])-2c_1(A), \quad\,
A \in \Gamma.
\]
It turns $CF (M;H)$ into a graded $\Lambda_\Gamma$-vector space.
For an admissible almost complex structure $J \in \cj$ define
the moduli space $\cm \left( [x,\bar{x}],[y,\bar{y}] \right)$
as the set of solutions $u$ of the Floer equation~\eqref{prob:floer}
for which
$\bar{x} \# u \# \bar{y}$ represents the trivial class in $\Gamma$. 
For generic choice of $J$
this moduli space is a smooth manifold of dimension 
$\mu \left( [x,\bar{x}] \right) - \mu \left( [y,\bar{y}] \right)$. 
Using Corollary~\ref{maximum}, the semi-positivity assumption
and the Floer--Gromov's compactness theorem 
one can prove that for generic choice of $J$
the moduli spaces $\cm([x,\bar{x}],[y,\bar{y}])$ for 
$\mu([x,\bar{x}])-\mu([y,\bar{y}])=1$ are compact, 
see \cite{HS}. 
We can thus set
\[
n \left( [x,\bar{x}], [y,\bar{y}] \right) \,:=\, 
\# \cm \left( [x,\bar{x}],[y,\bar{y}] \right)
\mod 2 .
\]
Define the Floer boundary operator $\pp_k \colon CF_k(M;H) \ra
CF_{k-1}(M;H)$ as the linear extension of
\[
\pp_k \left( [x,\bar{x}] \right) \, = 
\sum_{ \substack{ [y,\bar{y}] \in \widetilde{\cp}_H\\
\mu \left( [y,\bar{y}] \right) = k-1}}
n \left( [x,\bar{x}],[y,\bar{y}] \right) \, [y,\bar{y}]
\]
where $[x,\bar{x}] \in \widetilde{\cp}_H$ and 
$\mu([x,\bar{x}])=k$. Using again convexity, semi-positivity and Floer--Gromov
compactness one can prove that the right-hand side lies in
$CF_{k-1}(M;H)$, i.e., satisfies the required finiteness conditions. 
The boundary operator satisfies $\pp^2=0$, and its homology does not
depend on the regular pair $\left( H,J \right)$.
The resulting graded homology $HF_*(M)$ is a module over the
Novikov ring $\Lambda_\Gamma$.
Proceeding as in Section~\ref{iso} one constructs the PSS isomorphism
\[
\Phi \colon HM_* \left( M, \Lambda_\Gamma \right) 
:= HM_*(M) \otimes_{\ZZ_2} \Lambda_\Gamma   
\,\ra\, HF_* \left( M  \right) .
\]

\m
We are now going to explain how the Schwarz norm $\gg$ can be defined on
the level of functions.
For $H \in \ch$ the action spectrum $\Sigma_H$ is the set
\[
\Sigma_H \,=\, \left\{ \widetilde{\ca}_H \left( [x,\bar{x}] \right) \mid
[x,\bar{x}] \in \widetilde{\cp}_H \right\} .
\]
For $H \in \ch_{\reg}$ we define $c(H) \in \Sigma_H$ as is \eqref{def:cH}.
One shows as in the weakly exact case that $c$ is continuous with
respect to the Hofer norm, and so we can define $c$ on $\ch_c(M)$.
In order to see that $c$ satisfies the triangle
inequality~\eqref{triangle}, we notice that
the pair of pants product still defines a ring structure on Floer homology.
This ring structure is isomorphic to the ring structure on quantum homology
given by the quantum cup product. 
It is proved in \cite[Proposition 8.1.4]{MS2} 
that the image under the PSS isomorphism $\Phi$
of a point at which an admissible Morse function attains its 
single maximum is still the identity element in 
Floer homology endowed with the pair of pants ring structure.
Using this one shows as in \cite[Section~4]{Sch} that $c$ indeed
satisfies the triangle inequality.


\begin{lemma}  \label{l:app:pos}
Consider a compact split-convex symplectic manifold $(M, \go)$ with
$c_1(\go) =0$.
Assume that the time-independent and $C^2$-small Hamiltonian 
$H \in \ch_c(M)$ attains its maximum only at one point $p$, and that $p$
is a nondegenerate critical point.
Then
\[
c(H) \,=\, \max H .
\]
\end{lemma}

\proof
Since $c_1$ vanishes on $\pi_2 \left( M \right)$ and $H$ is
$C^2$-small, the Conley--Zehnder indices of the critical points of $H$
agree with their Morse indices. 
Since $H$ attains its maximum only at $p$, the point $p$ 
is the only critical point of index $2n$ and hence must be the
image of the PSS isomorphism $\Phi$ applied to the single point at which
an admissible Morse function attains its maximum. 
It follows that $c(H) = \max H$.   
\proofend

We define the Schwarz norm $\gg \colon \ch_c (M) \ra \RR$ by 
\[
\gg \left( H \right) \,=\, c \left( H \right) + c \left( H^- \right) .
\]
We do not study the relation between $\gg (H)$ and $\gg (K)$ if $\gf_H =
\gf_K$.
We do notice, however, that the continuity of $c$ and the fact that
$\Sigma_H \subset \RR$ has measure zero, \cite[Lemma 2.2]{Oh0},
imply that $\gg (H) = \gg(K)$ if $K$ is a time reparametrization of $H$.

\begin{theorem}  \label{t:app:norm}
Consider a compact split-convex symplectic manifold with $c_1(\go) =0$.
The Schwarz norm $\gg$ on $\ch_c(M)$ has the following properties.
\begin{itemize}
\item[(S1)]
$\gg ( 0 ) =0$ and $\gg (H)>0$ if $\gf_H \neq \idd$;
\item[(S2)]
$\gg ( H \Diamond K ) \le \gg (H) + \gg (K)$;
\item[(S3)]
$\gg ( H_\gt ) = \gg (H)$ for all $\gt \in \Sympcc (M)$;
\item[(S4)]
$\gg (H) = \gg \left( H^- \right)$;
\item[(S5)]
$\gg ( H ) \le \left\| H \right\|$.
\end{itemize}
Moreover, $\gg (H) \le 2\, \gg (K)$ if $\gf_K$ displaces $\supp \gf_H$.
\end{theorem}

\proof
Properties (S2), (S3), (S4) and (S5) are derived
as in the proof of Theorem~\ref{t:norm}.
Since $\gg (H)$ does not depend on the time parametrization of $H$,
since $c$ is continuous and since $\Sigma_H$ has measure zero,
we can prove the last statement by arguing as in the proof of
Proposition~\ref{p:phi2psi}.
The identity $\gg (0) =0$ follows from (S2) and (S5) with $H=K=0$.
The nondegeneracy of $\gg$ follows from Lemma~\ref{l:app:pos} and the last
statement.
\proofend

Theorem~\ref{t:app:norm} is strong enough to obtain the following version of
Theorem~\ref{t:dense} for rational split-convex symplectic manifolds with $c_1 =0$.

\begin{theorem}  \label{t:app:dense}
Assume that $(M, \go)$ is a rational split-convex symplectic manifold
with $c_1(\go) =0$ and index of rationality $\hbar$, and consider a
displaceable hypersurface $S$ of $M$ with displacement energy 
$e \left( S,M \right) < \hbar /2$.
Then for every thickening $\left( S_\eps \right)$ of $S$ and  
every $\gd > 0$ there exists $\eps \in \left[ -\gd, \gd \right]$
such that $\cp^\circ \left( S_\eps \right) \neq \emptyset$. 
\end{theorem}

\proof
We choose $K \in \ch_c(M)$ such that $\left\| K \right\| < \hbar /2$,
choose $i$ so large that $\supp \gf_K \subset M_i$, and define $H$ (with
$E=\hbar$) as in the proof of Theorem~\ref{t:dense}.
Properties (S1) and (S5) and the last statement in 
Theorem~\ref{t:app:norm} yield
\[
0 \,<\, \gg_{M_i} (H) \,\le\, 2\, \left\| K \right\| \,<\, \hbar .
\]
Since $H$ does not depend on time, it is obvious that $\Sigma_{H^-} = -
\Sigma_H$.
Notice that the contribution of a constant orbit $x \in \cp_H$ to
$\Sigma_H$ is $\hbar \ZZ$.
Together with $\gg_{M_i}(H) = c \left( H \right) + c \left( -H \right)$ we
conclude as in the proof of Theorem~\ref{t:dense} that 
$\cp^\circ \left( S_\eps \right) \neq \emptyset$ for some $\eps \in
\;]-\gd, \gd[$. 
\proofend


\begin{thebibliography}{99}



\bibitem{BC}
P.\ Biran and K.\ Cieliebak.
Lagrangian embeddings into subcritical Stein manifolds. 
{\it Israel J. Math.} {\bf 127} (2002) 221--244. 


\bibitem{BPS} 
P.\ Biran, L.\ Polterovich and D.\ Salamon.
Propagation in Hamiltonian dynamics and relative symplectic homology.
math.SG/0108134


\bibitem{Ci} 
K.\ Cieliebak.
Subcritical Stein manifolds are split.
math.DG/0204351







\bibitem{CIPP1} 
G.\ Contreras, R.\ Iturriaga, G.\ P.\ Paternain and M.\ Paternain.
Lagrangian graphs, minimizing measures and Ma\~n\'e's critical values. 
{\it Geom. Funct. Anal.} {\bf 8} (1998) 788--809. 


\bibitem{CIPP2} 
G.\ Contreras, R.\ Iturriaga, G.\ P.\ Paternain and M.\ Paternain.
The Palais-Smale condition and Ma\~n\'e's critical values. 
{\it Ann. Henri Poincaré} {\bf 1} (2000) 655--684. 


\bibitem{CMP} 
G.\ Contreras, L.\ Macarini and G.\ P.\ Paternain.
Periodic orbits for magnetic flows on surfaces.
Preprint 2002.
http://www.dpmms.cam.ac.uk/~gpp24/clor2.pdf



\bibitem{E} 
P.\ Eberlein.
Geometry of $2$-step nilpotent groups with a left invariant metric. 
{\it Ann. Sci. \'Ecole Norm. Sup.} {\bf 27} (1994) 611--660. 


\bibitem{E1}
Ya.\ Eliashberg. 
Topological characterization of Stein manifolds of dimension $>2$. 
{\it Internat. J. Math.} {\bf 1} (1990) 29--46.


\bibitem{E2}
Ya.\ Eliashberg. 
Symplectic geometry of plurisubharmonic functions. 
With notes by Miguel Abreu. NATO Adv. Sci. Inst. Ser. C
Math. Phys. Sci., 488, Gauge theory and symplectic geometry (Montreal,
1995), 49-67, Kluwer Acad. Publ., Dordrecht, 1997. 



\bibitem{EG}
Ya.\ Eliashberg and M.\ Gromov. Convex symplectic manifolds. In {\it
Several Complex Variables and Complex Geometry, Proceedings, Summer
Research Institute, Santa Cruz, 1989, Part 2}. Ed. by E.\ Bedford et
al.. Proc. Sympos. Pure Math. {\bf 52}, Amer. Math. Soc., Providence,
1991, 135-162.


\bibitem{F2}
A.\ Floer.
A relative Morse index for the symplectic action.
{\it Comm. Pure Appl. Math.} {\bf 41} (1988) 393--407.


\bibitem{F3}
A.\ Floer.
The unregularized gradient flow of the symplectic action.
{\it Comm. Pure Appl. Math.} {\bf 41} (1988) 775--813.

\bibitem{F1} 
A.\ Floer.    
Morse theory for Lagrangian intersections. 
{\it J. Diff. Geom.} {\bf 28} (1988) 513--547. 

\bibitem{F4}
A.\ Floer.
Witten's complex and infinite-dimensional Morse theory.
{\it J. Diff. Geom.} {\bf 30} (1989) 207--221.

\bibitem{F5}
A.\ Floer.
Symplectic fixed points and holomorphic spheres.
{\it Comm. Math. Phys.} {\bf 120} (1989) 575--611.


\bibitem{FHS} 
A.\ Floer, H.\ Hofer and D.\ Salamon.
Transversality in elliptic Morse theory for the symplectic action. 
{\it Duke Math. J.} {\bf 80} (1995) 251--292. 


\bibitem{F} 
U.\ Frauenfelder.
Floer homology of symplectic quotients and the Arnold--Givental
conjecture.
Diss. ETH No. 14981. Z\"urich 2003.


\bibitem{FS} 
U.\ Frauenfelder and F.\ Schlenk.
Slow entropy and symplectomorphisms of cotangent bundles.
Preprint ETH Z\"urich 2003.


\bibitem{FS2} 
U.\ Frauenfelder and F.\ Schlenk.
Spectral metrics and Hofer's geometry.
In preparation.





 


\bibitem{GT} 
D.\ Gilbarg and N.\ Trudinger.
{\it Elliptic partial differential equations of second order.} 
Second edition. Grundlehren der Mathematischen Wissenschaften {\bf 224}.
Springer-Verlag, Berlin, 1983. 



\bibitem{Gi0} 
V.\ Ginzburg.
On closed trajectories of a charge in a magnetic field. 
An application of symplectic geometry.
Contact and symplectic geometry (Cambridge, 1994), 131--148, 
{\it Publ. Newton Inst.} {\bf 8}, 
Cambridge Univ. Press, Cambridge, 1996. 


        
\bibitem{Gi}
V.\ Ginzburg.
Hamiltonian dynamical systems without periodic orbits. 
Northern California Symplectic Geometry Seminar, 35--48, 
{\it Amer. Math. Soc. Transl. Ser. 2}, {\bf 196}, 
Amer. Math. Soc., Providence, RI, 1999. 



\bibitem{GG}
V.\ Ginzburg and B.\ G\"urel.
Relative Hofer--Zehnder capacity and periodic orbits in twisted 
cotangent bundles.
math.DG/0301073. 


\bibitem{GK}
V.\ Ginzburg and E.\ Kerman.
Periodic orbits in magnetic fields in dimensions greater than two. 
Geometry and topology in dynamics 
(Winston-Salem, NC, 1998/San Antonio, TX, 1999), 113--121, 
{\it Contemp. Math.} {\bf  246}, Amer. Math. Soc., Providence, RI, 1999. 


\bibitem{GK2}
V.\ Ginzburg and E.\ Kerman.
Periodic orbits of Hamiltonian flows near symplectic extrema. 
{\it Pacific J. Math.} {\bf 206} (2002) 69--91.



\bibitem{GW}
D.\ Gromoll and J.\ Wolf.
Some relations between the metric structure and the algebraic structure
of the fundamental group in manifolds of nonpositive curvature. 
{\it Bull. Amer. Math. Soc.} {\bf 77} (1971) 545--552. 


\bibitem{Gr} 
M.\ Gromov. Pseudo-holomorphic curves in symplectic manifolds. 
{\it Invent. math.} {\bf 82} (1985) 307--347.


\bibitem{He}
D.\ Hermann.
Holomorphic curves and Hamiltonian systems in an open set with
restricted contact-type boundary.  
{\it Duke Math. J.} {\bf 103} (2000) 335--374.


\bibitem{H}
H.\ Hofer.
Estimates for the energy of a symplectic map. 
{\it Comment. Math. Helv.} {\bf 68} (1993) 48--72. 
 

\bibitem{HS}
H.\ Hofer and D.\ Salamon.
Floer homology and Novikov rings.
The Floer memorial volume, 483--524, 
{\it Progr. Math.} {\bf 133}.
Birkh\"auser-Verlag, Basel, 1995. 


\bibitem{HV}
H.\ Hofer and C.\ Viterbo.
The Weinstein conjecture in cotangent bundles and related results. 
{\it Ann. Scuola Norm. Sup. Pisa Cl. Sci. IV} {\bf 15} (1988) 411--445. 



\bibitem{HZ1} 
H.\ Hofer and E.\ Zehnder. 
Periodic solutions on hypersurfaces and a result by C. Viterbo. 
{\it Invent. Math.} {\bf 90} (1987) 1--9.


\bibitem{HZ2}
H.\ Hofer and E.\ Zehnder. 
A new capacity for symplectic manifolds. 
{\it Analysis, et cetera}, 405--427, 
Academic Press, Boston, MA, 1990. 


\bibitem{HZ} 
H.\ Hofer and E.\ Zehnder. 
{\it Symplectic Invariants and Hamiltonian Dynamics}. 
Birkh\"auser, Basel, 1994.


\bibitem{Ko} 
V.\ V.\ Kozlov. 
Calculus of variations in the large and classical mechanics. 
{\it Russian Math. Surveys} {\bf 40} (2) (1985) 37--71.


\bibitem{LS}
F.\ Laudenbach and J.-C.\ Sikorav. 
Hamiltonian disjunction and limits of Lagrangian submanifolds. 
{\it Internat. Math. Res. Notices} 1994, 161--168. 


\bibitem{L}
E.\ Lima. 
Orientability of smooth hypersurfaces and the Jordan-Brouwer separation
theorem. 
{\it Exposition. Math.} {\bf 5} (1987) 283--286. 



\bibitem{LM}
F.\ Lalonde and D.\ McDuff.
The geometry of symplectic energy. 
{\it Ann. of Math.} {\bf 141} (1995) 349--371. 


\bibitem{LP}
F.\ Lalonde and L.\ Polterovich.
Symplectic diffeomorphisms as isometries of Hofer's norm. 
{\it Topology} {\bf 36} (1997) 711--727.




\bibitem{Lu}
G.\ Lu
Periodic motion of a charge on a manifold in the magnetic fields.
math.DG/9905146.


\bibitem{Mac}
L.\ Macarini.
Hofer--Zehnder capacity and Hamiltonian circle actions.
math.SG/0205030. 
        

\bibitem{M1}
D.\ Mc\;\!Duff.
Symplectic manifolds with contact type boundaries. 
{\it Invent. Math.} {\bf 103} (1991) 651--671. 



\bibitem{MS} D.\ Mc\;\!Duff and D.\ Salamon. 
{\it Introduction to Symplectic Topology}. 
Oxford Mathematical Monographs, Clarendon Press, 1995.


\bibitem{MS2}
D.\ Mc\;\!Duff and D.\ Salamon. 
{\it $J$-holomorphic curves and quantum cohomology}. 
University Lecture Series {\bf 6}. 
American Mathematical Society, Providence, RI, 1994.



\bibitem{No}
S.\ P.\ Novikov.
The Hamiltonian formalism and a many-valued analogue of Morse theory. 
{\it Russian Math. Surveys} {\bf 37} (5) (1982) 1--56.

\bibitem{Oh0} 
Y.-G.\ Oh.
Chain level Floer theory and Hofer's geometry of the Hamiltonian
diffeomorphism group.
math.SG/0104243. 


\bibitem{Oh} 
Y.-G.\ Oh.
Normalization of the Hamiltonian and the action spectrum.
math.SG/0206090.


\bibitem{Oh2} 
Y.-G.\ Oh.
Mini-max theory, spectral invariants and geometry of the Hamiltonian
diffeomorphism group.
math.SG/0207214.



\bibitem{PP} 
G.\ Paternain and M.\ Paternain.
Critical values of autonomous Lagrangian systems. 
{\it Comment. Math. Helv.} {\bf 72} (1997) 481--499. 


\bibitem{PPS}
G.\ Paternain, L.\ Polterovich and K. Siburg.
Boundary rigidity for Lagrangian submanifolds, 
non-removable intersections, and Aubry-Mather theory.
math.SG/0207140. 



\bibitem{PSS} 
S.\ Piunikhin, D.\ Salamon and M.\ Schwarz.
Symplectic Floer-Donaldson theory and quantum cohomology. 
{\it Contact and symplectic geometry (Cambridge, 1994)}, 171--200, 
Publ. Newton Inst. {\bf 8}, Cambridge Univ. Press, Cambridge, 1996. 



\bibitem{P2}
L.\ Polterovich.
An obstacle to non-Lagrangian intersections. 
The Floer memorial volume, 575--586, 
Progr. Math. {\bf 133}, Birkh\"auser, Basel, 1995. 


\bibitem{P0}
L.\ Polterovich.
Geometry on the group of Hamiltonian diffeomorphisms. 
Proceedings of the International Congress of Mathematicians, 
Vol.\ II (Berlin, 1998). {\it Doc. Math.} {\bf 1998},
Extra Vol.\ II, 401--410.

\bibitem{P}
L.\ Polterovich.
{\it The geometry of the group of symplectic diffeomorphisms}.
Lectures in Mathematics ETH Z\"urich.
Birkh\"auser Verlag, Basel, 2001. 


\bibitem{P1} 
L.\ Polterovich. 
Growth of maps, distortion in groups and symplectic geometry.
{\it Invent. math.} {\bf 150} (2002) 655-686.


\bibitem{Sa} 
D.\ Salamon. 
Lectures on Floer homology. 
Symplectic geometry and topology (Park City, UT, 1997), 143--229, 
IAS/Park City Math. Ser. {\bf 7}, Amer. Math. Soc., Providence, RI, 1999. 


\bibitem{SZ} 
D.\ Salamon and E.\ Zehnder. 
Morse theory for periodic solutions of Hamiltonian systems and the
Maslov index. 
{\it Comm. Pure Appl. Math.} {\bf 45} (1992) 1303--1360. 


\bibitem{Sch1} M.\ Schwarz. 
{\it Morse homology}. Progress in Mathematics {\bf 111}.
Birkh\"auser Verlag, Basel, 1993.


\bibitem{Sch} M.\ Schwarz. 
On the action spectrum for closed symplectically
aspherical manifolds.
{\it Pacific J. Math.} {\bf 193} (2000) 419--461.


\bibitem{V1} C.\ Viterbo. 
Symplectic topology as the geometry of generating functions.
{\it Math. Ann.} {\bf 292} (1992) 685--710.


\bibitem{V2} C.\ Viterbo. 
Functors and computations in Floer homology with applications. I. 
{\it Geom. Funct. Anal.} {\bf 9} (1999) 985--1033.


\bibitem{Y}
S.\ Yau. 
On the fundamental group of compact manifolds of non-positive curvature.
{\it Ann. of Math.} {\bf 93} (1971) 579--585. 


\end{thebibliography}
\enddocument